\documentclass[11pt,reqno,tbtags]{amsart}
\usepackage{amssymb}
\numberwithin{equation}{section}

\newtheorem*{property*}{Property \csname @currentlabel\endcsname}

\newtheorem{theorem}{Theorem}[section]
\newtheorem{lemma}[theorem]{Lemma}
\newtheorem{proposition}[theorem]{Proposition}
\newtheorem{corollary}[theorem]{Corollary}

\theoremstyle{definition}

\newtheorem{remark}[theorem]{Remark}

\newtheorem*{ack}{Acknowledgement}

\theoremstyle{remark}

\newenvironment{romenumerate}{\begin{enumerate}
 }{\end{enumerate}}

\newcounter{oldenumi}
{\setcounter{oldenumi}{\value{enumi}}
\begin{romenumerate} \setcounter{enumi}{\value{oldenumi}}}
{\end{romenumerate}}

\newcounter{thmenumerate}

\newcounter{xenumerate}   

\newcommand{\refT}[1]{Theorem~\ref{#1}}
\newcommand{\refC}[1]{Corollary~\ref{#1}}
\newcommand{\refL}[1]{Lemma~\ref{#1}}
\newcommand{\refR}[1]{Remark~\ref{#1}}
\newcommand{\refS}[1]{Section~\ref{#1}}
\newcommand{\refP}[1]{Proposition~\ref{#1}}
\newcommand{\refand}[2]{\ref{#1} and~\ref{#2}}

\newcommand\marginal[1]{\marginpar{\raggedright\parindent=0pt\tiny #1}}

\begingroup
  \divide\count255 by 60
  \multiply\count255 by -60
  \advance\count255 by \time
  \ifnum \count255 < 10 \xdef\klockan{\the\count1.0\the\count255}
  \else\xdef\klockan{\the\count1.\the\count255}\fi
\endgroup

\newcommand\set[1]{\ensuremath{\{#1\}}}

\newcommand\bigpar[1]{\bigl(#1\bigr)}
\newcommand\Bigpar[1]{\Bigl(#1\Bigr)}

\newcommand\lrpar[1]{\left(#1\right)}

\newcommand\lrsqpar[1]{\left[#1\right]}
\newcommand\bigabs[1]{\bigl|#1\bigr|}

\newcommand\lrabs[1]{\left|#1\right|}
\def\rompar(#1){\textup(#1\textup)}    
\newcommand\xfrac[2]{#1/#2}

\newcommand\parfrac[2]{\Bigpar{\frac{#1}{#2}}}

\def\xexp(#1){e^{#1}}
\newcommand\ceil[1]{\lceil#1\rceil}
\newcommand\floor[1]{\lfloor#1\rfloor}

\newcommand\ntoo{\ensuremath{{n\to\infty}}}

\newcommand\bmin{\wedge}
\newcommand\bmax{\vee}

\newcommand\ie{i.e.\spacefactor=1000}
\newcommand\eg{e.g.\spacefactor=1000}

\newcommand\cf{cf.\spacefactor=1000}
\newcommand{\as}{a.s.\spacefactor=1000}

\newcommand\ii{\mathrm{i}}

\newcommand{\tend}{\longrightarrow}
\newcommand\dto{\overset{\mathrm{d}}{\tend}}
\newcommand\pto{\overset{\mathrm{p}}{\tend}}

\newcommand\eqd{\overset{\mathrm{d}}{=}}

\newcommand\bbR{\mathbb R}

\newcommand\bbZ{\mathbb Z}

\newcounter{CC} 
\newcommand{\CC}{\stepcounter{CC}\CCx} 
\newcommand{\CCx}{C_{\arabic{CC}}}     
\newcounter{cc}

\newcommand\E{\operatorname{\mathbb E{}}}
\renewcommand\P{\operatorname{\mathbb P{}}}
\newcommand\Var{\operatorname{Var}}
\newcommand\Cov{\operatorname{Cov}}

\newcommand\Bi{\operatorname{Bi}}
\newcommand\Be{\operatorname{Be}}

\newcommand\ga{\alpha}
\newcommand\gb{\beta}

\newcommand\gD{\Delta}

\newcommand\gam{\gamma}

\newcommand\gs{\sigma}
\newcommand\gss{\sigma^2}

\newcommand\cA{\mathcal A}

\newcommand\cM{\mathcal M}

\newcommand\ett[1]{\boldsymbol1[#1]} 

\def\[#1]{[\![#1]\!]}

\newcommand\qq{^{1/2}}
\newcommand\qqw{^{-1/2}}
\newcommand\qqq{^{1/3}}
\newcommand\qqqw{^{-1/3}}
\newcommand\qqa{^{2/3}}
\newcommand\qw{^{-1}}
\newcommand\qww{^{-2}}

\renewcommand{\=}{:=}

\newcommand\intot{\int_0^t}
\newcommand\intoo{\int_0^\infty}
\newcommand\intoooo{\int_{-\infty}^\infty}
\newcommand\oi{[0,1]}
\newcommand\oiq{[0,1)}
\newcommand\ooo{[0,\infty)}
\newcommand\oooo{\ensuremath{(-\infty,\infty)}}

\newcommand\dd{\,\textup{d}}

\newcommand\lhs{left hand side}
\newcommand\rhs{right hand side}
\newcommand\half{\frac12}
\newcommand\thalf{\tfrac12}
\newcommand\upto{\uparrow}
\newcommand\downto{\downarrow}

\newcommand{\afH}{af H\"allstr\"om}
\newcommand{\afHH}{\afH{} \cite{afH}}
\newcommand\Ai{\mathrm{Ai}}
\newcommand{\xnx}[1]{\ensuremath{X_{n,#1}}}
\newcommand{\xnm}{\ensuremath{\xnx{m}}}
\newcommand{\xnnt}{\ensuremath{\xnx{\floor{nt}}}}
\newcommand{\xnnh}{\ensuremath{\xnx{\ceil{n/2}}}}
\newcommand{\xx}{X^*}
\newcommand{\xxn}{\xx_n}
\newcommand{\xxx}{\bar X}
\newcommand{\xxxn}{{\xxx}_n}
\newcommand{\ynx}[1]{\ensuremath{Y_{n,#1}}}
\newcommand{\ynm}{\ensuremath{\ynx{m}}}
\newcommand{\ynnt}{\ensuremath{\ynx{\floor{2nt}}}}
\newcommand{\yy}{Y^*}
\newcommand{\yyn}{\yy_n}
\newcommand{\inx}[1]{I_{n,#1}}
\newcommand{\inm}{\inx{m}}
\newcommand{\sumkn}{\sum_{k=1}^n}
\newcommand{\sumkni}{\sum_{k=1}^{n-1}}
\newcommand\qn{}
\newcommand{\inq}{I\qn}
\newcommand{\inxx}[2]{\inq(#1;#2)}
\newcommand{\intx}[1]{\inxx{t}{#1}}
\newcommand{\intk}{\intx k}
\newcommand{\iinq}{I\qn'}
\newcommand{\iinxx}[2]{\iinq(#1;#2)}
\newcommand{\iintx}[1]{\iinxx{t}{#1}}
\newcommand{\iintk}{\iintx k}
\newcommand{\tinq}{{\tilde I}\qn}
\newcommand{\tinxx}[2]{\tinq(#1;#2)}
\newcommand{\tintx}[1]{\tinxx{t}{#1}}
\newcommand{\tintk}{\tintx k}
\newcommand{\tiinq}{{\tilde I}\qn'}
\newcommand{\tiinxx}[2]{\tiinq(#1;#2)}
\newcommand{\tiintx}[1]{\tiinxx{t}{#1}}
\newcommand{\tiintk}{\tiintx k}
\newcommand{\hinq}{{\widehat I}\qn}
\newcommand{\hinxx}[2]{\hinq(#1;#2)}
\newcommand{\hintx}[1]{\hinxx{t}{#1}}
\newcommand{\hintk}{\hintx k}
\newcommand{\hintkx}[1]{\hinxx{T_k}{#1}}
\newcommand{\hinsx}[1]{\hinxx{s}{#1}}
\newcommand{\hinsk}{\hinsx k}
\newcommand{\nnt}{\ensuremath{N_n(t)}}
\newcommand{\xn}{\ensuremath{X_n}}
\newcommand{\xnt}{\ensuremath{\xn(t)}}
\newcommand{\xnh}{\ensuremath{\xn(\thalf)}}
\newcommand{\yn}{\ensuremath{Y_n}}
\newcommand{\ynt}{\ensuremath{\yn(t)}}
\newcommand{\tx}[1]{\ensuremath{T_{(n;#1)}}}
\newcommand{\txm}{\tx m}
\newcommand{\tnt}{\tx{\floor{nt}}}
\newcommand{\tnh}{\tx{\ceil{n/2}}}
\newcommand{\xntnh}{\ensuremath{\xn\bigpar{\tnh}}}
\newcommand\gnp{\ensuremath{G(n,p)}}
\newcommand\gnm{\ensuremath{G(n,m)}}
\newcommand\ettn{\ensuremath{1,\dots,n}}
\newcommand{\sn}[1]{S_{n,#1}}
\newcommand{\sni}{\sn i}

\newcommand{\snx}[1]{\sn{#1}^*}
\newcommand{\hsn}[1]{\widehat S_{n,#1}}
\newcommand{\hsni}{\hsn i}
\newcommand{\tsn}[1]{\widetilde S_{n,#1}}
\newcommand{\snn}[1]{S_{n;#1}}
\newcommand{\snna}{\snn{\ga}}
\newcommand{\hsnn}[1]{\widehat S_{n;#1}}
\newcommand{\wn}[1]{W_{n,#1}}
\newcommand{\wni}{\wn i}
\newcommand{\wnj}{\wn j}
\newcommand\qv[1]{[#1,#1]}
\newcommand\qvv[2]{[#1,#2]}
\newcommand{\hZ}{\widehat Z}
\newcommand{\hZx}[1]{\widehat Z_{#1}}
\newcommand{\gsij}{\Sigma_{ij}}
\newcommand{\moo}{\cM_{\infty}}
\newcommand{\moox}[1]{\cM_{\infty #1}}
\newcommand{\qi}{^{(1)}}
\newcommand{\xnnqi}{X_n^{(1)*}}
\newcommand{\qd}{^{(d)}}

\newcommand{\tT}{\tilde T}
\newcommand{\oti}{{0\le t\le 1}}
\newcommand{\mnx}[1]{M_{n,#1}}
\newcommand{\mnxx}[1]{M_{n;#1}}
\newcommand{\mnj}{\mnx j}
\newcommand{\mn}{M_n}
\newcommand{\yoo}{\YX_\infty}
\newcommand{\YX}{V}
\newcommand{\xxxxn}{\xxx^*_n}
\newcommand{\xzxn}{X^{(1)*}_n}
\newcommand{\zz}{_*}
\newcommand{\ttt}{\tfrac13}

\newcommand{\Holder}{H\"older}

\newcommand{\maple}{\texttt{Maple}}

\newcommand\REM[1]{{\raggedright\texttt{[#1]}\par\marginal{XXX}}}

\newcommand\urladdrx[1]{{\urladdr{\def~{{\tiny$\sim$}}#1}}}

\begin{document}
\title[Sorting by subintervals and the maximum number of runs]
{Sorting using complete subintervals and the maximum number of runs in
a randomly evolving sequence.}

\date{January 10, 2006}

\author{Svante Janson}
\address{Department of Mathematics, Uppsala University, PO Box 480,
SE-751~06 Uppsala, Sweden}
\email{svante.janson@math.uu.se}
\urladdrx{http://www.math.uu.se/~svante/}


\begin{abstract} 
We study the space requirements of a sorting algorithm where only
items that at the end will be adjacent are kept together.
This is equivalent to the following combinatorial problem:
Consider a string of fixed length $n$ that starts as a
string of 0's, and then evolves by changing each 0 to 1, with the
$n$ changes done in random order. What is the maximal number of runs of 1's?

We give asymptotic results for the distribution and mean. It turns out
that, as in many problems involving a maximum, the maximum is
asymptotically normal, with fluctuations of order $n\qq$, and to the
first order well approximated by the number of runs at the instance
when the expectation is maximized, in this case when half the elements
have changed to 1; there is also a second order term of order $n\qqq$.

We also treat some variations, including priority queues.

The proofs use methods originally developed for random graphs.
\end{abstract}

\maketitle

\section{Introduction}\label{S:intro}

Gunnar \afHH{} considered, as indicated at the end of his paper, 
the following algorithm for sorting an unordered pile of student exams in
alphabetic order. 
(It is said that he used this procedure himself.)

The exams are taken one by one from the input.
The first exam is put in a new pile.
For each following exam ($x$, say), if the name on it 
is immediately preceding the name on an exam $y$ at the top of one
of the piles, the new exam $x$ is put on top of $y$. (The professor
knows the names of all the students, and can thus see that there are
no names between $x$ and $y$.)
Similarly, if the name on $x$ is immediately succeeding
the name on an exam $z$ at the bottom of a pile, $x$ is put under
$z$. If both cases apply, with $y$ on top of one pile and $z$ at the
bottom of another, the two piles are merged with $x$ inserted between
$z$ and $y$. Finally, if there is no pile matching $x$ in one of these
ways, $x$ is put in a new pile.

The algorithm thus maintains a list of sorted piles, each being an
interval without gaps of the set of exams. At the end, there is a
single sorted pile.

The problem is the space requirement of this algorithm; more
precisely, the maximum number of sorted piles during the execution.
The input is assumed to be in random order, so this is a random
variable, and we are interested in its mean and distribution.

\begin{remark}
  As a sorting method, this algorithm has drawbacks. First, it
  requires that all names are known from the beginning; mathematically
  it can be seen as sorting the numbers \ettn.
Secondly, the space requirement turns out to be quite high, see
  below. This also implies that the number of comparisons necessary
  for each insertion is high, of the order of $n$.
The algorithm might be useful when blocks of sorted items can be
  manipulated as easily as individual items, and we do not want to
  make insertions inside the blocks, for example when sorting physical
  objects that are to be glued together in order.
\end{remark}

\afHH{} gave the following mathematical reformulation, where we also
introduce some notation.
Consider a deck of $n$ cards numbered \ettn{} in random order,
and a sequence of $n$ places with the same numbers in order.
Take the cards one by one and put them at their respective places.
When we have placed $m$ cards, $0\le m\le n$, we see $\xnm$
``islands'', \ie{} uninterrupted blocks of cards.
What is $\xxn\=\max_m\xnm$?

Alternatively,
we can  use the language of parking cars, which is popular
for some related problems in computer science:
$n$ cars park, one by one, on $n$ available places along a street;
each car parks at a random free place. What is the maximum
number of  uninterrupted blocks of cars during the process?

Let, for $n\ge1$, $0\le m\le n$ and $1\le k\le n$,
the indicator $\inm(k)$ be 1 if the item (exam or card) with number
$k$ is one of the
$m$ first in the input, and 0 otherwise.
Thus, $\xnm$ is the number of runs of 1's in the random sequence
$\inm(1),\dots,\inm(n)$ of $n-m$ 0's and $m$ 1's.
We can express \xnm{} algebraically as
\begin{equation}
  \label{a2}
  \begin{split}
\xnm
&=\inm(1)+\sumkni(1-\inm(k))\inm(k+1)
\\
&=m-\sumkni \inm(k)\inm(k+1).	
  \end{split}
\end{equation}

If the input is given by the permutation $\gs$ of \set\ettn, so that
item $k$ has position $\gs\qw(k)$,
\begin{equation*}
  \inm(k)=\ett{\gs\qw(k)\le m},
\end{equation*}
where $\ett{\dots}$ denotes the indicator of the indicated event.
We assume that $\gs$ is a (uniformly chosen) random permutation; thus
so is $\gs\qw$.
Hence, each random sequence $(\inx{m}(k))_{k=1}^n$ is uniformly
distributed over all $\binom nm$ possibilities; moreover,
for each $m<n$ we obtain $(\inx{m+1}(k))_{k=1}^n$ from $(\inx{m}(k))_{k=1}^n$
by changing a single randomly chosen 0 to 1, this random choice being
uniform among the $n-m$ 0's, and independent of the previous history.

It is easy to see that $\E \xnm=m(n-m+1)/n$, see \eqref{e1}; it follows
that the maximum of $\E\xnm$ for a given $n$ is attained  for
$m=\ceil{n/2}$, and that $\E\xnx{\ceil{n/2}}>n/4$.
Since obviously
\begin{equation}
  \label{af1}
\E\xxn=\E\max_m\xnm\ge\max_m\E\xnm,
\end{equation}
this yields $\E\xxn>n/4$ as observed by
\afHH. Moreover, he observed that $\E\xxn$ is subadditive, and thus the
limit
\begin{equation*}
  \gamma\=\lim_{\ntoo} \E\xxn/n
\end{equation*}
exists and equals $\inf_{n} \E\xxn/n$; he further showed that
$1/4\le\gamma\le1/3$, where the lower 
bound comes from \eqref{af1}.
Based on simulations with $n=13$ and $n=52$, \afHH{} concluded that $\gamma$
seems to be very close to or equal to 1/4.
We will show that, indeed, $\gam=1/4$.
We also show that the distribution of $\xxn$ is asymptotically normal,
with a variance of order $n$.

\begin{theorem}
  \label{T1}
As \ntoo,
\begin{equation}\label{t1}
  n\qqw\bigpar{\xxn-n/4} \dto N(0,1/16),
\end{equation}
with convergence of all moments.
In particular,
\begin{align*}
  \E\xxn&=n/4+o(n\qq),
\\
\Var\xxn&=n/16+o(n).
\end{align*}
\end{theorem}

This theorem says that to the first order, the maximum number of piles
(runs) $\xxn$ behaves like the number $\xnm$ with $m=\ceil{n/2}$. A
more refined analysis shows that the difference $\xxn-\xnx{\ceil{n/2}}$
is of order $n\qqq$. 
Let $B(t)$, $-\infty<t<\infty$, be a standard two-sided Brownian
motion; thus $B(0)=0$ and $B(t)$, $t\ge0$, and $B(-t)$, $t\ge0$, are
two independent 
Brownian motions.
\begin{theorem}
  \label{T2}
As \ntoo,
\begin{equation}\label{t2}
  n\qqqw\bigpar{\xxn-\xnnh}\dto 
\thalf \YX,
\end{equation}
where the random variable $\YX$ is defined by 
$\YX\=\max_{t} \bigpar{B(t)-t^2/2}$,
and
\begin{equation*}
  \E\xxn=
\E\xnx{\ceil{n/2}}+ \thalf\E \YX n\qqq +o(n\qqq)
=\tfrac14 n+ \thalf\E \YX n\qqq +o(n\qqq).
\end{equation*}
\end{theorem}
The random variable $\YX$ is studied by Barbour \cite{B75}, Daniels
and Skyrme \cite{DS85} and Groeneboom \cite{Groeneboom}.
Note that $0<\YX<\infty$ a.s.
We have, see \cite{DS85} (using \maple{} to improve the numerical
values in \cite{B75,B81,DS85,D89}), with $\Ai$ the Airy function,
\begin{equation*}
 \E \YX=-\frac{2\qqqw}{2\pi}\intoooo \frac{\ii y\dd y}{\Ai(\ii y)^2} 
\approx 0.996193.
\end{equation*}

The numerical values $\xx_{13}\approx4.22$ and $\xx_{52}\approx14.66$
found experimentally by af H\"allstr\"om \cite{afH} differ from $n/4$ by
about 18\% and 10\% less than the correction term $\thalf\E \YX n\qqq$
in \refT{T2}, which is a reasonable agreement for such rather small $n$.

\begin{remark}\label{Rcyclic}
  af H\"allstr\"om \cite{afH} considered also the cyclic case, when we
  regard \set{1,\dots,n} as a circle, which sometimes is slightly
  simpler to study
  because of the greater symmetry.
In this case we define $\inm(k)$ for all $k\in\bbZ$ by
  $\inm(k+n)\=\inm(k)$, \ie{} we
  interpret $k$ modulo $n$, and we sum to $n$ in \eqref{a2}.
Since the number of runs in the linear and cyclic version differ by at
  most 1, all our asymptotic results remain the same, 
and we will only consider the linear case.
(Moreover, the cyclic case with $n$ items corresponds exactly to the
  linear with $n-1$ by fixing the last element, see \cite{afH}.)
\end{remark}

We prove these theorem by studying asymptotics of the entire (random) process
$(\xnm)_{m=0}^n$. The natural time here is $m/n$, 
so we take $m=\floor{nt}$ for $0\le t\le 1$
and consider the process $\xnnt$ with a continuous parameter
$t\in\oi$.
The following theorem shows that this process asymptotically is Gaussian.
(The space $D\oi$ is defined in \refS{Ssn}, see \cite{Bill} for a
detailed treatment.)

\begin{theorem}
  \label{TXm}
As \ntoo, in the space $D\oi$ of functions on $\oi$,
\begin{equation}
  \label{txm1}
n\qqw\bigpar{\xnnt-nt(1-t)}\dto Z(t),
\end{equation}
where $Z$ is a continuous Gaussian
process on $\oi$ 
with mean $\E Z(t)=0$ and covariances
\begin{align}
\label{txm2}
  \E\bigpar{Z(s)Z(t)}&=s^2(1-t)^2,&&
0\le s\le t\le 1.
\end{align}
\end{theorem}

The behaviour of $\xxn$ shown in Theorems \refand{T1}{T2}, with an
asymptotic normal distribution with a mean of order $n$ and random
fluctuations of order $n\qq$, and with a second order term for the
mean of order $n\qqq$, is common for this type of random variables
defined as the maximum of some randomly evolving process. 
For various examples, both combinatorial and others, and general results
see for example Daniels \cite{D74,D89}, Daniels and Skyrme \cite{DS85},
Barbour \cite{B75,B81} and Louchard, Kenyon and Schott \cite{LKS97}.
Indeed, paraphrasing the explanations in these papers,
in many such problems, the
first order asymptotic of a random process $X_n(t)$ 
(after suitable scaling) is a
deterministic function $f(t)$, say, defined on a compact interval $I$
(typically scaled to be $\oi$ as here).
Hence the first order asymptotic of the maximum of the process is just
the maximum of this function $f$. Moreover, it is often natural to
expect that the random fluctuations around this function $f(t)$
asymptotically form a Gaussian process $G(t)$; this is then a second
order term of smaller order as in our \refT{TXm}. 
If we assume that $f$ is continuous on
$I$ and has a unique maximum at a point $t_0\in I$, then the maximum
of the process $X_n(t)$ is attained close to $t_0$, so the first order
approximation of the maximum is the constant $f(t_0)=\max_t f(t)$, while
the next approximation is just $X_n(t_0)$, giving a normal limit law
as in our \refT{T1}. The Gaussian fluctuations in this limit have mean
0, so in order to find the next term for the mean $\E\xxn$, we 
study more closely the difference $\max_t X_n(t)-X_n(t_0)$ by studying
the difference $X_n(t)-X_n(t_0)$ close to $t_0$. Assuming that $t_0$
is an interior point of $I$ and that $f$ is
twice differentiable at $t_0$ with $f''(t_0)\neq0$, we can locally at
$t_0$ approximate $f$ by a parabola and $G(t)-G(t_0)$ by a two-sided
Brownian motion (with some scaling), and thus $\max_t X_n(t)-X_n(t_0)$
is approximated by a scaling constant times the variable $\YX$ above,
see Barbour \cite{B75} and, in our case, \refC{CX2} below.
In the typical case where the mean of $X_n(t)$ is of order $n$ and the
Gaussian fluctuations are of order $n\qq$, it is easily seen that
the correct scaling gives, as in \refT{T2} above, a correction to
$\E\xxn$ of order $n\qqq$, see  \cite{B75,D74,D89} and \refS{Sproofs}.

  The method used in the present paper is a simple adaption of the
  method used in \cite{SJ79} and \cite{SJ94} to study the number of
   subgraphs of a given isomorphism type in a random graph.
   These papers study the random 
  graphs $\gnp$ and $\gnm$ that can be constructed by random deletion
  of edges in the complete graph $K_n$ (with the deletions being
  independent for $\gnp$ and such that a fixed number of edges are
  deleted for \gnm). The method applies more generally to random
  graphs constructed  by random edge deletions in these ways from any
  fixed initial graph $F_n$. The problem treated in this paper can be
  regarded
as an instance of this when the initial graph is the path $P_n$ with
  $n$ edges.
In particular, \refT{T1} corresponds to \cite[Theorem 24]{SJ94}, which
gives the asymptotic distribution of the maximum number of induced
  subgraphs of a given type during the evolution of \gnp{} or \gnm;
see also \cite[Theorem 33]{SJ94} (isolated edges) and 
\cite[Theorem 17]{SJ94} (a general result) for related results. 
Conversely, 
we expect that these
results for random graphs can be complemented by the analogues of
\refT{T2} above, using the method of proof in the present paper,
but we have not verified the details.

Our method applies also to other problems. 
First, let $\xnm\qi$ be the number of piles with a single exam (runs
with a single 1) in the process studied above. Then we obtain similar
results for the maximum $\xnnqi\=\max_m\xnm\qi$, see \refS{Sfurther}.
The same applies to the number $\xnm\qd$ of piles with any other fixed
number $d$ of exams (runs of a fixed length $d$).

Another example is given by \emph{priority queues}, where Louchard
\cite{L87} and Louchard,  Kenyon and Schott \cite{LKS97} have proved
asymptotic results very similar to the Theorems \ref{T1}--\ref{TXm}
above. 
In particular, they found the same asymptotic covariance \eqref{txm2}
except for a normalizing constant.
(See also 
Flajolet, Fran\c{c}on and Vuillemin \cite{FFV80}
and Flajolet, Puech and Vuillemin \cite{FPV86}
for combinatorial results on generating functions involving Hermite
polynomials; these results, however, do not easily yield asymptotics.)

Priority queues can be defined as follows. Suppose that $n$ items are
to be temporarily stored (or processed); let item $i$ arrive at time
$A_i$ and be deleted at time $D_i$.
We assume that the $2n$ times $A_i$ and $D_i$ are distinct; thus they
can be arranged in a sequence of the $2n$ events $A_i$ and $D_i$, with
$A_i$ coming before $D_i$ for each $i$. We assume further, as our
probabilistic model, that all $(2n)!/2^n$ such sequences are equally probable.
Ignoring the labels, we can equivalently consider sequences of $n$ $A$
and $n$ $D$ (or $n$ $+$ and $n$ $-$), where each $A$ is paired with a
$D$ coming later; there is a 1--1 correspondence between such
sequences and pairings of $1,\dots,2n$ into $n$ pairs, and there are
$(2n-1)!!=(2n)!/(2^nn!)$ such sequences (with pairings), again taken
with equal probability.

Let, for $m=0,\dots,2n$, $\ynm$ be the number of items stored after
$m$ of these events, \ie{} the number of $A$'s minus the
number of $D$'s among the $m$ first events, and let $\yyn\=\max_{0\le
m\le2n}\ynm$. 
The sequence $(\ynm)_0^{2n}$ is a Dyck path, but note that its
distribution is not uniform; for a given Dyck path (or a given
sequence of $A$ and $D$ without labels), the number of ways to pair
a given $D$ with a preceding $A$, \ie{} the number of ways to choose
which item to delete, equals the current number of items stored before
this deletion. Thus, the weight of the Dyck path equals the product of
these numbers $\prod_{m:\ynx{m+1}<\ynm}\ynm$.
Alternatively, which better explains the name priority queue,
we can keep the stored items in a list showing the
order in which they eventually will be deleted; then there is only one
choice for each deletion but each new item can be inserted in $Y+1$
ways if there are $Y$ items stored before the insertion, and thus
$Y+1$ after it; hence the weight can also be written as
$\prod_{m:\ynm>\ynx{m-1}}\ynm$. (It is easily to see directly that the
two products are equal.)

We will in \refS{Spriority} show how our method applies to priority
queues, and explain why we obtain the same asymptotic results as for
$\xnm$ and $\xxn$. 
(Note that there is no exact correspondence for
finite $n$, since the natural sample spaces have $n!$ elements for
$\xnm$ but $(2n-1)!!$ elements for $\ynm$.) Again, we can regard the
problem as an instance of subgraph counts for randomly deleting edges
from a given initial graph $F_n$; in this case taking $F_n$ to
be a multigraph consisting of $n$ double edges.

A third example
is a model suggested by Van Wyk and Vitter \cite{VV86} as a model for
hashing with lazy deletion, and further studied by
Louchard \cite{L88} and Louchard,  Kenyon and Schott \cite{LKS97}.
In this model, $n$ item arrives and are deleted as above, but now the
arrival and deletion times $A_i$ and $D_i$ are random numbers,
with the $n$ pairs $(A_i,D_i)$ mutually independent and each pair distributed
as $(T_i\bmin \tT_i,T_i\bmax \tT_i)$, where $T_i$ and $\tT_i$ are
independent random variables uniformly distributed on \oi.
(We use $\bmin$ and $\bmax$ as notations for $\min$ and $\max$ of two numbers.)
We let $\ynt$ be the number of items present at time $t$, and again
we are especially interested in its maximum $\max_t\ynt$. 
Again, the
asymptotic 
results for the maximum found by Louchard,  Kenyon and Schott \cite{LKS97} 
are the same as in our
Theorems \ref{T1} and \ref{T2}, except for a constant factor, while
the asymptotic result for the process $\yn(t)$ found by Louchard
\cite{L88} differs somewhat from the one in \refT{TXm}; it
corresponds instead to the one in \refC{CX} below.
Indeed, as explained by Kenyon and Vitter \cite{KV91},
see also \refS{Spriority},
this model can be seen as a priority queue with randomized
times for insertions and deletions, which explains why the results for
the maximum are the same as for priority queues.

We assume in the sequel that $n\ge2$, to avoid some trivialities.
All unspecified limits are as \ntoo.
We use the standard notations $\pto$ and $\dto$ for convergence in
probability and distribution, respectively, of random variables, and
\as{} for \emph{almost surely}, \ie{} with probability 1.

\section{Randomizing time}\label{Srand}

We will use the standard method of randomizing the time.
More precisely, we let $T_1,\dots,T_n$ be independent random
variables, each uniformly distributed on $(0,1)$.
We interpret $T_k$ as the time item $k$ arrives, and note that \as{}
there are no ties.
We define
\begin{equation*}
 \intk=\ett{T_k\le t},  
\end{equation*}
\ie, $\intk=1$ if item $k$ has arrived by time
$t$. 
We further define $\nnt$ as the number of items that have arrived at time $t$,
and $\xnt$ as the number of runs of 1's at time $t$, \ie, \cf{} \eqref{a2},
\begin{align}
\label{mt}
\nnt&=\sumkn\intk,
\\
\xnt
&=\intx{1}+\sumkni\bigpar{1-\intk}\intx{k+1}
  \label{xnt1}
\\&
=\nnt-\sumkni \intk\intx{k+1}.
  \label{xnt}
\end{align}

Clearly, the items arrive in random order, so the process remains the
same except that the insertions occur at the random times
$\tx1,\dots,\tx{n}$, where \tx{j} is the $j$:th order statistic of
$T_1,\dots,T_n$. We thus have $\inm(k)=\inxx{\tx{m}}{k}$ and
$\xnt=\xnm$ when $\tx{m}\le t<\tx{m+1}$ (with $\tx0\=0$ and
$\tx{n+1}\=1$ for convenience).
In particular,
\begin{equation}\label{julie}
  \xxn=\max_{0\le t\le1}\xnt.
\end{equation}
Note that $X_n(0)=\xnx0=0$ and $X_n(1)=\xnx{n}=1$.

The importance of this randomization is that the variables $\intk$,
$k=1,\dots,n$, are independent (both for a fixed $t$ and as
stochastic processes, \ie{} as random functions of $t$).
For every $n$, $k$ and $t\in\oi$,
\begin{equation}
  \label{a4}
\P\bigpar{\intk=1}
=
\P(T_k\le t)=t,
\end{equation}
\ie{} $\intk$ has the Bernoulli distribution $\Be(t)$.
\xnt{} thus is the number of runs of 1 in a sequence of \emph{independent}
0's and 1's, each with the distribution $\Be(t)$.
Furthermore, the number of items sorted at time $t$ is $\nnt\sim\Bi(n,t)$.

Define further, for $0\le t\le 1$, the centralized variables
\begin{equation}  \label{ii}
  \iintk\=\intk-\E\intk=\intk-t
\end{equation}
and the sums
\begin{align}  \label{sn1}
  \sn1(t)&\=\sumkn\iintk=\nnt-\E\nnt=\nnt-nt,
\\
\label{sn2}
  \sn2(t)&\=\sumkni\iintk\iintx{k+1}.
\end{align}
Thus $\sn1(0)=\sn2(0)=\sn1(1)=\sn2(1)=0$ and $\E\sn1(t)=\E\sn2(t)=0$
for all $t\in\oi$.
We have
\begin{align}
  \nnt \label{nnt}
&=
\sumkn\bigpar{\iintk+t}
=\sn1(t)+nt,
\\
\sumkni\intk\intx{k+1}
&=
\sumkni\bigpar{\iintk+t}\bigpar{\iintx{k+1}+t} \notag
\\
&=\sn2(t)+t\bigpar{2\sn1(t)-\iintx1-\iintx{n}}+(n-1)t^2,
\notag
\end{align}
and thus from \eqref{xnt} the representation
\begin{equation}
  \label{sofie}
  \begin{split}
\xnt
&=n(t-t^2)+t^2+(1-2t)\sn1(t)-\sn2(t)+t\iintx1+t\iintx{n}
\\
&=nt(1-t)+(1-2t)\sn1(t)-\sn2(t)+R_n(t),	
  \end{split}
\raisetag{12pt}
\end{equation}
where $R_n(t)\=t^2+t\iintx1+t\iintx{n}$ and thus $|R_n(t)|\le3$.

We will in \refS{Ssn} study the asymptotic distribution of the
stochastic processes (\ie, random functions) $\sn1(t)$ and $\sn2(t)$;
our main results then follow easily from \eqref{julie} and \eqref{sofie}.

Note that for any fixed $t$, the variables $\iintk$ are independent
and have means 0; hence the terms in the sums in \eqref{sn1} and
\eqref{sn2} have means and all covariances 0. (They are thus
orthogonal in $L^2$.) It follows immediately that
\begin{align}
  \Var\lrpar{\sn1(t)} 
&= n\E\bigpar{\iintx1}^2
= n\Var\bigpar{\intx1}
=nt(1-t),
\label{csn11}
\\
 \Var\lrpar{\sn2(t)} 
&= (n-1)\E\bigpar{\iintx1\iintx2}^2
= (n-1)\bigpar{\Var\bigpar{\intx1}}^2
\notag\\&
=(n-1)t^2(1-t)^2,
\label{csn22}
\\
\Cov\bigpar{\sn1(t),{}&\sn2(t)}
=0.
\label{csn12}
\end{align}

\section{Exact results} \label{Sexact}

We first give some exact results for finite $n$.
It is easy to find the exact distribution of $\xnm$ for given $n$ and
$m$, see for example Stevens \cite{Stevens} or Mood \cite{Mood}.
For $m\ge1$ and $k\ge1$
we have $\xnm=k$ if there are $k$ runs of 1's separated by
$k-1$ runs of 0's and possibly preceded and/or succeeded by additional
runs of 0's. Considering the bivariate generating function for such
sequences of arbitrary length, we easily find 
\begin{align*}
  \P(\xnm=k)&=
[x^my^{n-m}] \parfrac{x}{1-x}^k \parfrac{y}{1-y}^{k-1}\parfrac{1}{1-y}^2
\\
&=[x^{m-k}y^{n-m-k+1}](1-x)^{-k}(1-y)^{-k-1}
\\
&=\binom{m-1}{k-1}\binom{n-m+1}{k}.
\end{align*}
The mean can be computed from this \cite{Mood}, \cite{afH}, but simpler from
\eqref{a2}:
\begin{equation}
  \label{e1}
  \begin{split}
\E\xnm
&=m-\sumkni\E\bigpar{ \inm(k)\inm(k+1)}
\\
&=m-(n-1)\frac{m(m-1)}{n(n-1)}=\frac{m(n-m+1)}{n}.
  \end{split}
\end{equation}
A similar computation of the variance yields, omitting the details,
\begin{align*}
  \Var \xnm = \frac{m(m-1)(n-m)(n-m+1)}{n^2(n-1)}.
\end{align*}

If we instead randomize the insertion times as in \refS{Srand} and
consider the process at a fixed time $t$, we have by \eqref{xnt},
\eqref{a4} and the independence of $\intk$ for $k=1,\dots,n$,
\begin{equation}\label{e2}
  \E \xnt=t+\sumkni(1-t)t= nt(1-t)+t^2.
\end{equation}
Similarly, using \eqref{xnt1}, again omitting details,
\begin{align}\label{e3}
  \Var \xnt
=nt(1-t)(1-3t+3t^2)+t^2(1-t)(3-5t).
\end{align}

To find the exact distribution of $\xxn$ seems much more complicated.
Exact values of $\P(\xxn=h)$ are easily calculated for small $n$, see
\afHH, but we do not know any general formula. It would be interesting
to find such a formula by combinatorial methods.

\section{The asymptotic distribution of $\sn1(t)$ and $\sn2(t)$}
\label{Ssn}

To state our results on the asymptotic distribution of the
stochastic processes $\sn1(t)$ and $\sn2(t)$, we need a suitable
topological space of functions. We use, for an interval
$I\subseteq\bbR$, 
the standard space $D(I)$ of right-continuous functions on $I$ 
that have left-hand 
limits, equipped with the Skorohod topology. For a precise definition
of this (metrizable) topology, see \eg{}
Billingsley \cite{Bill} ($I=\oi$),
Jacod and Shiryaev \cite{JS} ($I=\ooo$), 
Kallenberg \cite[Appendix A.2]{Kallenberg} ($I=\ooo$), 
or Janson \cite{SJ94}. 
For our purposes it is sufficient to know that 
if $f$ is continuous on $I$, then $f_n\to f$ in $D(I)$
if and only if $f_n\to f$ uniformly on every compact subinterval.
In particular, if $I$ is compact, for example $I=\oi$, 
and $f$ is continuous on $I$, then $f_n\to f$ in $D(I)$
if and only if $f_n\to f$ uniformly.

Our main result on the asymptotic global behaviour of $\sn1(t)$ and
$\sn2(t)$ then can be stated as follows. 
\begin{theorem}
  \label{TS}
As \ntoo, in $D\oi$,
\begin{align}
  \label{c1}
n\qqw\sn1(t)&\dto Z_1(t),
\\
  \label{c2}
n\qqw\sn2(t)&\dto Z_2(t),
\end{align}
jointly, where $Z_1$ and $Z_2$ are two independent continuous Gaussian
processes on $\oi$ with means $\E Z_1(t)=\E Z_2(t)=0$ and covariances
\begin{align}
  \E\bigpar{Z_1(s)Z_1(t)}&=s(1-t),&&
 0\le s\le t\le 1,
\label{c3}
\\
  \E\bigpar{Z_2(s)Z_2(t)}&=s^2(1-t)^2,&&
0\le s\le t\le 1.
\label{c4}
\end{align}
\end{theorem}
Thus, $Z_1$ is a standard Brownian bridge, and
the limit \eqref{c1} is just the well-known theorem 
that the
empirical distribution function asymptotically is distributed as a
Brownian bridge, see \eg{} Billingsley \cite[Theorem 16.4]{Bill}.

The proof of \refT{TS}, and of all other results in this section, are
postponed to \refS{Sproofs}. 

Using \eqref{sofie}, \refT{TS} yields the asymptotic distribution of
the process $\xnt$. 
\begin{corollary}
  \label{CX}
As \ntoo, in $D\oi$,
\begin{equation}
  \label{cx}
n\qqw\bigpar{\xnt-nt(1-t)}\dto Z(t),
\end{equation}
where $Z$ is a continuous Gaussian
process on $\oi$ with mean $\E Z(t)=0$ and covariances,
for $ 0\le s\le t\le 1$,
\begin{align}
  \E\bigpar{Z(s)Z(t)}
&=s(1-2s)(1-t)(1-2t)+s^2(1-t)^2
  \label{cx1}
\\
&=s(1-t)(1-s-2t+3st).
\label{cx2}
\end{align}
\end{corollary}

In particular, this implies the limit \eqref{cx} for each fixed $t$,
with $\Var(Z(t))=t(1-t)(1-3t+3t^2)$, which also follows more easily
from \eqref{xnt1}, \eqref{e2}, \eqref{e3} and the Central Limit Theorem
for 1-dependent sequences.

These results are stated using the randomized insertions described in
\refS{Srand}. 
We can also return to the original deterministic insertion times and
obtain asymptotics of the discrete process $(\xnm)_{m=0}^n$, which
yields \refT{TXm} stated in the introduction. Note that the limit
processes in \refT{TXm} and \refC{CX} are different, due to the additional
random variation introduced when randomizing the time.
(The variance of the limit in \refT{TXm} is strictly smaller
than in \refC{CX} at every $t\notin\set{0,\half,1}$.)

We will also need a moment estimate. It is easy to see that 
$n\qq\sni(t)$ has moments that are bounded as \ntoo, for every fixed
$t\in\oi$. We extend that to the supremum over all $t$.

\begin{theorem}
  \label{TSx}
Let $\snx{i}\=\sup_{0\le t\le1} |\sni(t)|$
for $i=1,2$. 
Then, for each
fixed $r>0$,
$\E\bigpar{\snx{i}}^r = O(n^{r/2})$.
\end{theorem}

We are primarily interested in the maximum $\xxn$ of $\xnt$. It is
evident from \refC{CX} that the maximum is attained close to the
maximum point of $t(1-t)$, \ie, close to $t=1/2$. We use a magnifying
glass and study the processes close to $t=1/2$ in greater detail. The
correct scaling turns out to be $t=\half+xn\qqqw$, and we have the
following asymptotic behaviour on that scale.

\begin{theorem}
  \label{TS2}
As \ntoo, in $D(-\infty,\infty)$,
\begin{align}
  \label{w1}
n\qqqw\bigpar{\sn1(\thalf+xn\qqqw)-\sn1(\thalf)}&\dto B_1(x),
\\
  \label{w2}
n\qqqw\bigpar{\sn2(\thalf+xn\qqqw)-\sn2(\thalf)}&\dto2\qqw B_2(x),
\end{align}
jointly, where $B_1$ and $B_2$ are two independent Brownian motions on
\oooo.
Furthermore, for any fixed $A<\infty$ and $i=1,2$,
\begin{equation}
  \E \max_{|x|\le A} \lrpar{\sni(\thalf+xn\qqqw)-\sni(\thalf)}^2 
= O(n^{2/3}).
  \label{w3}
\end{equation}
\end{theorem}

\begin{corollary}
  \label{CX2}
As \ntoo, in $D(-\infty,\infty)$,
\begin{align}
  \label{erika}
n\qqqw\bigpar{X_n(\thalf+xn\qqqw)-X_n(\thalf)}\dto 2\qqw B(x)-x^2,
\end{align}
where $B$ is a Brownian motion on
\oooo.
\end{corollary}

\section{time-reversal}\label{Stime}
In the proofs below, we will introduce factors that blow up at the
endpoint $t=1$. To see that there is no real problem at this endpoint,
we will use a time reversal trick which enables us to transfer results
from the other endpoint.

If we replace each $T_k$ by $1-T_k$, 
which of course has the same distribution, then 
$\intk$ becomes $1-\inxx{1-t}k$, except at the jump point, and thus,
see \eqref{ii}--\eqref{sn2},
$\iintk$ becomes $-\iinxx{1-t}k$ and 
$\sni(t)$ becomes $(-1)^i\sni(1-t)$, again excepting the jump points.
To be precise, let for a
function $f$ on $\oi$, $f(t-):=\lim_{s\upto t}f(s)$ (when this exists),
with $f(0-)\=f(0)$. Then $\sni(t)$ becomes $(-1)^i\sni((1-t)-)$ under
this time-reversal, and thus
\begin{equation}\label{time}
  \sni(t)\eqd(-1)^i\sni((1-t)-),
\end{equation}
as functions in $D\oi$ and jointly for $i=1,2$.

\section{Proofs}\label{Sproofs}

The proofs are based on martingale theory, in particular a
continuous time martingale limit theorem by Jacod and Shiryaev \cite{JS}.   
We will use the \emph{quadratic variation} $[X,X]_t$
of a martingale $X$ (in continuous time) and its bilinear
extension $[X,Y]_t$ to two martingales $X$ and $Y$. 
For a general definition see \eg{} \cite{JS};
for us it will suffice to know that,
if $X$ and $Y$ are martingales of pathwise finite variation, then 
\begin{equation}\label{x1}
  [X,Y]_t=\sum_{0<s\le t} \gD X(s) \gD Y(s),
\end{equation}
where $\gD X(s)\=X(s)-X(s-)$ is the jump of $X$ at $s$
and, similarly, $\gD Y(s)\=Y(s)-Y(s-)$.
The sum in \eqref{x1} is formally uncountable, but in reality
countable since there is only a countable number of jumps; in the
applications below, the sum will be finite. 

A real-valued martingale $X(s)$ on $[0,t]$ is an $L^2$-martingale if
and only if  
$\E [X,X]_t<\infty$ and $\E|X(0)|^2<\infty$,
and then
\begin{equation}
  \label{x21}
  \E |X(t)|^2 = \E[X,X]_t+\E |X(0)|^2.
\end{equation}

We will use the following general result based on \cite{JS}; 
see \cite[Proposition 9.1]{SJ154} for a detailed proof
(for $I=\ooo$; the general case is the same).
(See also \cite{SJ79} and
\cite{SJ94} for similar versions).

\begin{proposition}\label{P:JS}
Let $I=[a,b]$ or $I=[a,b)$, with $-\infty<a<b\le\infty$.
Assume that for each $n$, 
$\cM_n(t)=(\cM_{ni}(t))^q_{i=1}$ is a 
$q$-dimensional martingale on $I$
with $\cM_n(a)=0$, and that 
$\Sigma(t)=(\gsij(t))_{i,j=1}^q$ 
is a  \rompar(non-random) 
continuous matrix-valued function on $I$
such that
for every fixed $t\in I$ 
and $i,j=1,\dots,q$,
\begin{align}
&[\cM_{ni}, \cM_{nj}]_t\pto \gsij(t)
\quad\text{as \ntoo,}
\label{js1}
\\
&\sup_n\E [\cM_{ni},\cM_{ni}]_t <\infty .
\label{js2}
\end{align}
Then $\cM_n\dto \moo$ as $n\to\infty$, in 
$D(I)$, where
$\moo=(\moox i)_{i=1}^q$ 
is a continuous $q$-dimensional Gaussian martingale with $\E \moo(t)=0$ and
covariances
\begin{equation*}
\E \bigpar{\moox i(s)\moox j(t)}=\gsij (s\bmin t),\qquad s,t\in I.
\end{equation*}
In other words, the components $\cM_{ni}(t)$ converge jointly to $\moox
i(t)$ in $D(I)$.
\end{proposition}

\begin{remark}
  By \eqref{x21}, \eqref{js2} is equivalent to
$\sup_n\E |\cM_n (t)|^2 <\infty$, the form used in \eg{} \cite{SJ154}.
\end{remark}

\begin{proof}[Proof of \refT{TS}]
We first construct martingales from $\sn1(t)$ and $\sn2(t)$.
We define, for $0\le t<1$,
\begin{align}
  \hintk&\=\frac{\iintk}{1-t}
=
\begin{cases}
  1, & \intk=1,
\\
-\xfrac{t}{(1-t)}, & \intk=0;
\end{cases}
\notag
\\
  \hsn1(t)&\=\sumkn\hintk=(1-t)\qw\sn1(t);
  \label{ems1}
\\
  \hsn2(t)&\=\sumkni\hintk\hintx{k+1}=(1-t)\qww\sn2(t).
  \label{ems2}
\end{align}
We have 
$\E\hintk=0$ and
\begin{equation}
  \label{emma}
\E\bigpar{\hintk^2}
=\Var\bigpar{\hintk}=(1-t)\qww\Var\bigpar{\intk}=\frac{t}{1-t}.
\end{equation}

It is easily checked that each $\hintk$ is a martingale on $[0,1)$
\cite[Lemma 2.1]{SJ94}; since these martingales for different $k$
are independent, the products $\hintk\hintx{k+1}$ are martingales too,
and thus $\hsn1(t)$ and $\hsn2(t)$ are martingales on $[0,1)$ with
  $\hsn1(0)=\hsn2(0)=0$. 
To calculate their quadratic variations and covariation, note that
$\gD\hintk=(1-t)\qw$ when $t=T_k$ and 0 otherwise.
Further, 
with $\hintx0\=\hintx{n+1}\=0$,
\begin{align*}
\gD \hsn1(t)&=\sumkn \gD\hintk,
\\
\gD \hsn2(t)&=\sumkn \gD\hintk\bigpar{\hintx{k-1}+\hintx{k+1}},
\end{align*}
and thus
\begin{align}
  \label{b11}
\qv{\hsn1}_t&=\sum_{s\le t}\sumkn\bigabs{\gD\hinsk}^2
=\sumkn\frac1{(1-T_k)^2}\ett{T_k\le t},
\\
  \label{b12}
\qvv{\hsn1}{\hsn2}_t
&=\sum_{s\le t}\sumkn\bigabs{\gD\hinsk}^2\bigpar{\hinsx{k-1}+\hinsx{k+1}}
\notag\\&
=\sumkn\frac1{(1-T_k)^2}\bigpar{\hintkx{k-1}+\hintkx{k+1}}\ett{T_k\le t},
\\
  \label{b22}
\qv{\hsn2}_t
&=\sum_{s\le t}\sumkn\bigabs{\gD\hinsk}^2\bigpar{\hinsx{k-1}+\hinsx{k+1}}^2
\notag\\&
=\sumkn\frac1{(1-T_k)^2}\bigpar{\hintkx{k-1}+\hintkx{k+1}}^2\ett{T_k\le t}.  
\end{align}

Hence, since the $T_k$ are independent and uniformly distributed on
\oi, and using \eqref{emma},
\begin{align}
  \label{be11}
\E\qv{\hsn1}_t
&= n\intot\frac{\dd s}{(1-s)^2}=n\lrsqpar{\frac1{1-s}}_0^t 
=n\frac{t}{1-t},
\\
  \label{be12}
\E\qvv{\hsn1}{\hsn2}_t
&=\sumkn
\intot\frac{\dd s}{(1-s)^2}\E\bigpar{\hinsx{k-1}+\hinsx{k+1}}
=0,
\\
  \label{be22}
\E\qv{\hsn2}_t
&=\sumkn
\intot\frac{\dd s}{(1-s)^2}\E\bigpar{\hinsx{k-1}+\hinsx{k+1}}^2
\notag\\
&=((n-2)\cdot2+2\cdot1)
\intot\frac{\dd s}{(1-s)^2}\frac{s}{1-s}
=(n-1)\frac{t^2}{(1-t)^2}
\end{align}
(Indeed, these formulas also follow directly from
\eqref{csn11}--\eqref{csn12} by \eqref{ems1}, \eqref{ems2} and \eqref{x21}
together with its polarized version for two martingales.)

Moreover, the $k$:th and $l$:th terms in the sums in
\eqref{b11}--\eqref{b22} are independent when $|k-l|>2$, and each term
is $O\bigpar{(1-t)\qww}$. Hence, for $i,j\in\set{1,2}$,
\begin{equation}
  \label{b7}
\Var\bigpar{\qvv{\hsn i}{\hsn j}_t}
= O\bigpar{n(1-t)^{-4}}.
\end{equation}
Define now, for $i=1,2$ and $0\le t<1$,
\begin{equation}\label{tsn}
  \tsn i(t)\=n\qqw\hsn i(t).
\end{equation}
By \eqref{be11}--\eqref{b7}, for every fixed $t\in\oiq$,
\begin{align*}
\qv{\tsn1}_t
&\pto\frac{t}{1-t},
\\
\qvv{\tsn1}{\tsn2}_t
&\pto0,
\\
\qv{\tsn2}_t
&\pto
\frac{t^2}{(1-t)^2}.
\end{align*}
\refP{P:JS} thus applies with $I=\oiq$, 
with \eqref{js2} verified by \eqref{tsn}, \eqref{be11} and \eqref{be22},
which shows that
\begin{equation}\label{hsz}
  n\qqw\hsn i(t)=  \tsn i(t)
\dto
\hZ_i(t),
\qquad i=1,2,
\end{equation}
jointly in $D\oiq$, where $\hZ_1(t)$ and $\hZ_2(t)$ are continuous
Gaussian processes on $\oiq$ with means 0 and covariances, for $0\le s\le t<1$,
\begin{align}
  \label{hz11}
\E \bigpar{\hZ_1(s)\hZ_1(t)}
&=\frac{s}{1-s},
\\
  \label{hz12}
\E \bigpar{\hZ_1(s)\hZ_2(t)}
&=0,
\\
  \label{hz22}
\E \bigpar{\hZ_2(s)\hZ_2(t)}
&=\frac{s^2}{(1-s)^2}.
\end{align}
Note that \eqref{hz12} implies that $\hZx1$ and $\hZx2$ are independent.

We define $Z_i(t)\=(1-t)^i \hZ_i(t)$ for $t\in\oiq$, and
$Z_i(1)\=0$. Then 
\eqref{hsz} implies, by \eqref{ems1} and \eqref{ems2},
that \eqref{c1} and \eqref{c2}
hold jointly in $D\oiq$. Furthermore, \eqref{hz11}--\eqref{hz22}
imply that $Z_1$ and $Z_2$ have the covariances \eqref{c3} and \eqref{c4}.

It remains to extend this from $\oiq$ to $\oi$.
We use the time-reversal trick in \refS{Stime} and have by
\eqref{time} and the result just shown
\begin{equation*}
n\qqw\sni((1-t)-) 
\eqd(-1)^i n\qqw \sni(t)
\dto (-1)^i Z_i(t)\eqd Z_i(t)
\end{equation*}
in $D\oiq$, and thus
$
n\qqw\sni(t)\dto Z_i(1-t)$
in $D(0,1]$. Clearly $Z_i(1-t)\eqd Z_i(t)$, as processes on $\oi$;
since $Z_i$ is continuous on $\oiq$, this implies continuity at 1 too,
and thus $Z_i$ is continuous on \oi.
We have shown that the limits \eqref{c1} and \eqref{c2} hold (jointly)
in both $D\oiq$ and $D(0,1]$, which easily implies 
that they hold in $D\oi$ too,
see \eg{} \cite[Lemma 2.3]{SJ94}.
\end{proof}

\begin{proof}[Proof of \refC{CX}]
  Immediate by \refT{TS} and \eqref{sofie}, 
with the limit $Z(t)\=(1-2t)Z_1(t)-Z_2(t)$,
since $n\qqw R_n(t)\to0$ uniformly.
The covariances \eqref{cx1}--\eqref{cx2} of $Z(t)$ follow from
\eqref{c3}, \eqref{c4} and the independence of $Z_1$ and $Z_2$.
\end{proof}

\begin{proof}[Proof of \refT{TSx}]
By \Holder's inequality, it suffices to prove the result when $r$ is
an even integer. 
Since 
$\snx i 
\le \sup_{0\le t\le1/2}\lrabs{\sni(t)}+\sup_{1/2\le
  t\le1}\lrabs{\sni(t)}$ and
the time-reversal symmetry \eqref{time} implies
\begin{equation}\label{dx2}
  \sup_{1/2\le t\le1}\lrabs{\sni(t)}
\eqd
  \sup_{0\le t\le1/2}\lrabs{\sni(t)},
\end{equation}
it is sufficient to consider
$\sup_{0\le t\le1/2}\lrabs{\sni(t)}$.
Moreover, $|\sni(t)|\le\bigabs{\hsni(t)}$, and by Doob's maximal
inequality for martingales, see \eg{} \cite[Proposition 7.16]{Kallenberg},
\begin{equation}
  \label{dx3}
\E\Bigpar{\sup_{0\le t\le1/2}\lrabs{\sni(t)}}^r
\le
\E\Bigpar{\sup_{0\le t\le1/2}\lrabs{\hsni(t)}}^r
\le C_r
\E\Bigpar{\lrabs{\hsni(\thalf)}}^r,
\end{equation}
for some constant $C_r$ ($=(r/(r-1))^r$).

Finally, $\hsn1(\thalf)$ is the sum of $n$ independent random
variables $\hinxx{\thalf}{k}$, each with values $\pm1$ and mean 0, and
it is easily verified that, with $r=2\ell$,
\begin{equation}
  \label{dx1}
\E\Bigpar{\lrabs{\hsni(\thalf)}}^r
=O(n^\ell)=O(n^{r/2}).
\end{equation}
Similarly, $\hsn2(\half)$ is the sum of the $n-1$ random variables 
$\hinxx{\half}{k}\hinxx{\half}{k+1}$; these variables too have values
$\pm1$ and mean 0; moreover, it is easily verified that they too are
independent. Hence $\hsn2(\thalf)\eqd S_{n-1,1}(\thalf)$, and
\eqref{dx1} implies the same estimate for $\hsn2(\thalf)$ too.

The result follows by this, \eqref{dx3} and \eqref{dx2}.
\end{proof}

\begin{proof}[Proof of \refT{T1}]
We claim that
\refC{CX} implies that 
\begin{equation}
  \label{d0}
n\qqw\Bigpar{\max_{0\le t\le1}X_n(t)-\tfrac14 n}
\dto Z(\thalf)\=-Z_2(\thalf),
\end{equation}
which gives \eqref{t1} by \eqref{c4}.
(We could use \refT{TXm} instead.)
The argument was sketched in the introduction, 
and this is an application of \cite[Theorem 16]{SJ94}, 
but for completeness we give the details in our case.
We may for simplicity use 
the Skorohod coupling theorem \cite[Theorem 4.30]{Kallenberg}, which
says that we can assume that \eqref{cx} holds with convergence \as{}
and not just in distribution.
Thus, 
for (almost) every point in our probability space, 
$n\qqw\bigpar{\xnt-nt(1-t)}\to Z(t)$ in $D\oi$, which since
$Z(t)$ is continuous means uniform convergence on \oi. 
In other words, uniformly in $t\in\oi$,
\begin{equation}\label{dd}
 \xnt=nt(1-t)+n\qq Z(t)+o(n\qq)
=\tfrac14n-n(\thalf-t)^2+n\qq Z(t)+o(n\qq).
\end{equation}
In particular,
\begin{equation}\label{d1}
  \xxn \ge \xn(\thalf)= \tfrac14 n+n\qq Z(\thalf)+o(n\qq).
\end{equation}
Conversely, \eqref{dd} yields for $|t-\thalf|<n^{-1/8}$, since $Z$ is
continuous, 
\begin{equation}\label{d2}
  \xnt\le \tfrac14n+n\qq Z(t)+o(n\qq)
=\tfrac14n+n\qq Z(\thalf)+o(n\qq),
\end{equation}
and for $|t-\thalf|\ge n^{-1/8}$, since $Z$ is bounded,
\begin{equation}\label{d3}
  \xnt
\le \tfrac14n-n^{1-1/4}+O(n\qq)
\le \tfrac14n+n\qq Z(\thalf)
\end{equation}
for large $n$. It follows from \eqref{d1}, \eqref{d2} and \eqref{d3}
that
$$\xxn=\tfrac14n+n\qq Z(\thalf)+o(n\qq),$$
and \eqref{d0} follows.

To prove moment convergence in \eqref{t1}, it is, as is well-known,
see \eg{} \cite[Theorems 5.4.2 and 5.5.9]{Gut}, enough to prove that
for each fixed $r>0$, the $r$:th absolute moment of the \lhs{} is
bounded, as \ntoo. By \eqref{sofie}, 
\begin{equation*}
  \begin{split}
\lrabs{\xxn-\tfrac14 n}
&=\bigabs{\sup_t \xnt-\sup_t nt(1-t)}
\le \sup_t \bigabs{\xnt-nt(1-t)} 
\\&
\le \snx1+\snx2+3,	
  \end{split}
\end{equation*}
and the required estimate follows by \refT{TSx}.
\end{proof}

\begin{proof}[Proof of \refT{TXm}]
Recall the order statistics \tx{m} from \refS{Srand}. Since
$\xnm=\xn(\tx m)$, we are studying the process 
$\xnnt=\xn(\tnt)$. The idea of the proof is to use the functional
limit results just shown and
replace $t$ by the random time \tnt.
Note first that $N_n(\tx m)=m$ and thus
\begin{equation}
  \label{p4a}
N_n(\tnt)=\floor{nt}=nt+O(1).
\end{equation}
By \eqref{sn1},
\begin{equation}\label{p4b}
  \begin{split}
  \sup_{0\le m\le n}\bigabs{N_n(\txm)/n-\txm}
&\le
  \sup_\oti\bigabs{N_n(t)/n-t}
=\sup_\oti\bigabs{n\qw\sn1(t)}
\\&
=n\qw\snx1,	
  \end{split}
\raisetag\baselineskip
\end{equation}
which by \eqref{c1} (or \refT{TSx}, or the 
Glivenko--Cantelli theorem \cite[Proposition 4.24]{Kallenberg})
tends to 0 in probability. 
Thus, by \eqref{p4a},
\begin{equation}\label{p4}
\sup_\oti\bigabs{t-\tnt}
\le
  \sup_\oti\bigabs{N_n(\tnt)/n-\tnt}+n\qw
\pto0.
\end{equation}
The proof of \refC{CX} shows that \eqref{cx} holds jointly with
\eqref{c1} and \eqref{c2}, with $Z(t)=(1-2t)Z_1(t)-Z_2(t)$.
Furthermore, by \eqref{sn1}, $N_n(t)/n=t+\sn1(t)/n$, and a Taylor
expansion of the function $t\mapsto nt(1-t)$ yields
\begin{equation*}
  N_n(t)\bigpar{1-N_n(t)/n}
=nt(1-t)+(1-2t)\sn1(t)-\sn1(t)^2/n.
\end{equation*}
Consequently, by \eqref{c1}, in $D\oi$,
still jointly with \eqref{cx}, 
\begin{equation*}
n\qqw\Bigpar{ N_n(t)\bigpar{1-N_n(t)/n}-nt(1-t)}
\dto(1-2t)Z_1(t),
\end{equation*}
and subtracting this from \eqref{cx}
yields
\begin{equation}\label{p5}
n\qqw\Bigpar{\xnt- N_n(t)\bigpar{1-N_n(t)/n})}
\dto Z(t)-(1-2t)Z_1(t)=-Z_2(t).
\end{equation}
Because \eqref{p4} holds and $Z_2(t)$ is continuous, we may 
replace $t$ by \tnt{} on the \lhs; 
for a formal verification of this we may again use
the Skorohod coupling theorem \cite[Theorem 4.30]{Kallenberg}
and thus assume that \eqref{p4} and  \eqref{p5} hold \as, \ie{} that
the functions in
\eqref{p4} and \eqref{p5} converge uniformly on \oi{} to their limits.
Consequently,
\begin{equation}\label{p6}
  n\qqw\Bigpar{\xnnt-N_n(\tnt)\bigpar{1-N_n(\tnt)/n}}
\dto -Z_2(t),
\end{equation}
which by \eqref{p4a} yields \eqref{txm1} with $Z(t)=-Z_2(t)$.
\end{proof}

The fact that the terms with $\sn1$ cancel in the proof above is no
coincidence. $\sn1$ measures by \eqref{sn1} the random fluctuations
introduced by used random insertion times $T_k$, and it is very
intuitive that this term will appear in the limits for $\xnt$ but
not for $\xnm$.
A theorem verifying that this cancellation happens in general in a
situation closely related to the one studied here is given in
\cite[Theorem 7]{SJ94}.

\begin{proof}[Proof of \refT{TS2}]
Fix $A>0$, and define for $n>(2A)^3$ and $x\in[0,2A]$,
\begin{equation}
  \label{w0}
\wni(x)\=\hsni\bigpar{\thalf+(x-A)n\qqqw}
-\hsni\bigpar{\thalf-An\qqqw}.
\end{equation}
Then $\wni$ is a martingale on $[0,2A]$
with $\wni(0)=0$,
and its quadratic variation is by \eqref{x1}
\begin{equation}
  \label{ww1}
\qv{\wni}_x=\qv{\hsni}_{\half+(x-A)n\qqqw}-\qv{\hsni}_{\half-An\qqqw}.
\end{equation}
Hence, by \eqref{be11}--\eqref{be22}, for $0\le x\le 2A$,
\begin{align}
\E\qv{\wn1}_x&=n\int_{\half-An\qqqw}^{\half+(x-A)n\qqqw} \frac{\dd s}{(1-s)^2}
=n\qqa x (4+O(n\qqqw)),
\label{ww11}
\\
\E\qvv{\wn1}{\wn2}_x&=0,
\label{ww12}
\\
\E\qv{\wn2}_x&=(2n-2)\int_{\half-An\qqqw}^{\half+(x-A)n\qqqw} 
 \frac{ s\dd s}{(1-s)^3}
=2n\qqa x (4+O(n\qqqw)).
\label{ww22}
\end{align}
Moreover, by \eqref{ww1} and \eqref{b7}, for $n>(4A)^3$, say,
\begin{equation*}
  \Var\qvv{\wni}{\wnj}_x=O(n).
\end{equation*}
Consequently, \refP{P:JS} applies to $n\qqqw\wni$, and shows that
in $D[0,2A]$ and jointly for $i=1,2$, 
\begin{equation}
  \label{ww2}
n\qqqw\wni(x)\dto W_i(x),
\end{equation}
where $W_1$ and $W_2$ are independent Gaussian stochastic processes
with means 0 and 
\begin{align*}
  \E\bigpar{W_1(x)W_1(y)}&=4x, &
  \E\bigpar{W_2(x)W_2(y)}&=8x, &
&0\le x\le y\le 2A.
\end{align*}
In other words, $W_1(x)=2B_1(x)$ and $W_2(x)=\sqrt8\,B_2(x)$, where
$B_1$ and $B_2$ are independent Brownian motions on $[0,2A]$. 
We may assume that $B_1$ and $B_2$ actually are independent two-sided Brownian
motions defined on the entire real line. Note that
$B_i(x+A)-B_i(A)\eqd B_i(x)$ (as processes on $\bbR$). Hence we can
make a translation and obtain from \eqref{w0} and \eqref{ww2}, 
in $D[-A,A]$ and jointly for $i=1,2$,
\begin{multline}\label{www}
n\qqqw\bigpar{\hsni\bigpar{\thalf+xn\qqqw}-\hsni\bigpar{\thalf}}
=n\qqqw\bigpar{\wni(x+A)-\wni(A)}
\\
\dto 2^{(i+1)/2}\bigpar{B_i(x+A)-B_i(A)}
\eqd	
2^{(i+1)/2}B_i(x).
\end{multline}

By \eqref{ems1} and \eqref{ems2} we further have,  
uniformly for $n>(4A)^3$ and $x\in[-A,A]$,
\begin{align}
\sni\bigpar{\thalf+xn\qqqw}-\sni\bigpar{\thalf}  
=
(\thalf-xn\qqqw)^i\hsni\bigpar{\thalf+xn\qqqw}-(\thalf)^i\hsni\bigpar{\thalf}
\notag
\\
=
2^{-i}\Bigpar{\hsni\bigpar{\thalf+xn\qqqw}-\hsni\bigpar{\thalf}}
+O(n\qqqw\snx i).
\label{wx}
\end{align}
and thus, using \eqref{www} and \refT{TSx}, in $D[-A,A]$ and jointly
for $i=1,2$,
\begin{equation*}
n\qqqw\bigpar{\sni\bigpar{\thalf+xn\qqqw}-\sni\bigpar{\thalf}}
\dto 2^{-i}
2^{(i+1)/2}B_i(x)
=
2^{(1-i)/2}B_i(x).
\end{equation*}  
Since convergence in $D[-A,A]$ for every $A>0$ implies convergence in
$D(-\infty,\infty)$, this proves \eqref{w1} and \eqref{w2}.

For the second moment estimate \eqref{w3}, we 
first note that \eqref{x21}, \eqref{ww11} and \eqref{ww22}  show that,
for each fixed $A$,
\begin{equation}\label{wx1}
 \E |\wni(2A)|^2 = \E\qv{\wni}_{2A}=O(n\qqa), 
\end{equation}
and thus by Doob's maximal inequality
\cite[Proposition 7.16]{Kallenberg},
\begin{equation*}
 \E\bigpar{ \max_{0\le x\le 2A}|\wni(x)|^2} = O(n\qqa).
\end{equation*}
Hence, by \eqref{w0} and translation again,
\begin{multline*}
\E\max_{|x|\le A}\bigabs{\hsni\bigpar{\thalf+xn\qqqw}-\hsni\bigpar{\thalf}}^2
=\E\max_{|x|\le A}\bigabs{\wni(x+A)-\wni(A)}^2
\\
\le
4\E\max_{0\le x\le 2A}\bigabs{\wni(x)}^2
=O(n\qqa),
\end{multline*}
and \eqref{w3} follows by \eqref{wx} and \refT{TSx}.
\end{proof}

\begin{remark}
  We can extend \eqref{w3} to arbitrary powers $r>0$, with the
  estimate $O(n^{r/3})$, by the argument above with the
  Burkholder--Davis--Gundy inequalities 
\cite[Theorem 26.12]{Kallenberg} replacing \eqref{wx1}; we omit the details.
\end{remark}

To study $\xnt$ close to $t=1/2$, we rewrite \eqref{sofie} as, 
for $|x|\le1/2$,
\begin{equation}
  \label{cecil}
\xn\bigpar{\thalf+x}
=\tfrac14n-nx^2-2x\sn1\bigpar{\thalf+x}-\sn2\bigpar{\thalf+x}
+R_n\bigpar{\thalf+x}.
\end{equation}
Hence, still for $|x|\le1/2$, 
\begin{multline}
  \label{cecilia}
\xn\bigpar{\thalf+x}-\xn\bigpar{\thalf}
=-nx^2-2x\sn1\bigpar{\thalf+x}
-\bigpar{\sn2\bigpar{\thalf+x}-\sn2\bigpar{\thalf}}
\\
+R_n\bigpar{\thalf+x}-R_n\bigpar{\thalf}
\end{multline}
and thus, for $|x|\le n\qqq/2$, recalling $|R_n(t)|\le3$,
\begin{multline}
  \label{cecilib}
n\qqqw\Bigpar{\xn\bigpar{\thalf+xn\qqqw}-\xn\bigpar{\thalf}}
=-x^2-2n^{-2/3}x\sn1\bigpar{\thalf+xn\qqqw}
\\
-n\qqqw\bigpar{\sn2\bigpar{\thalf+xn\qqqw}-\sn2\bigpar{\thalf}}
+O\bigpar{n\qqqw}.
\end{multline}

\begin{proof}[Proof of \refC{CX2}]
Fix $A>0$.
Then \eqref{erika} follows in $D[-A,A]$ by \eqref{cecilib}, \eqref{w2} and
\refT{TSx} (which implies $n^{-2/3}\snx i\pto0$),
with $B(x)\=-B_2(x)$.
Since $A>0$ is arbitrary, this yields convergence in
$D(-\infty,\infty)$.
\end{proof}

Let $x_+\=x\bmax0$.

\begin{lemma}
  \label{Lbound}
Let $x_1>0$ and
suppose that $\cM(x)$ is a martingale on $[0,x_1]$ 
with $\cM(0)=0$ such that
for some constant $K$ and all
$x\in[0,x_1]$
\begin{equation*}
  \Var \cM(x) \le Kx.
\end{equation*}
Then, for every $a>0$ and $x_0\in(0,x_1]$,
\begin{equation*}
  \E\Bigpar{\max_{x_0\le x\le x_1}\bigpar{\cM(x)-ax^2}_+} \le \frac{4K}{ax_0}.
\end{equation*}
\end{lemma}

\begin{proof}
We may for convenience extend $\cM$ to a martingale on $\ooo$ by
letting $\cM(x)\=\cM(x_1)$ for $x>x_1$.
Let $y>0$ and $t>0$. Then, by Kolmogorov-Doob's inequality
\cite[Proposition 7.16]{Kallenberg},
\cite[Theorem 10.9.1]{Gut},
\begin{align*}
  \P\Bigpar{\sup_{x\in[y,2y]}\bigpar{\cM(x)-ax^2}>t}
&\le
  \P\Bigpar{\sup_{x\in[0,2y]}{\cM(x)}>t+ay^2}
\\&
\le
\frac{\E \cM(2y)^2}{(t+ay^2)^2}
\le
\frac{2Ky}{(t+ay^2)^2}.
\end{align*}
Integrating with respect to $t$ from 0 to $\infty$ yields
\begin{align*}
  \E\bigpar{\sup_{x\in[y,2y]}\bigpar{\cM(x)-ax^2}_+}
\le
\intoo\frac{2Ky}{(t+ay^2)^2}
=\frac{2Ky}{ay^2}
=\frac{2K}{ay},
\end{align*}
and the result follows by summing over $y=2^k x_0$, $k=0,1,\dots$.
\end{proof}

\begin{proof}[Proof of \refT{T2}]
We begin by showing that we can replace $\xnnh$ by $\xnh$ in the
statement.
By \eqref{p4b} and \refT{TSx},
\begin{equation*}
 n\qqq\bigabs{N_n(\tnh)/n-\tnh}
\le n^{-2/3}\snx1 \pto0
\end{equation*}
and thus by \eqref{p4a}
\begin{equation}
  \label{q1}
n\qqq\bigabs{\tnh-\thalf}\pto0.
\end{equation}
It now follows from \refC{CX2}, arguing as for \eqref{p6} and
using \eqref{q1} and the fact that the limit in \eqref{erika} is continuous,
that we can
substitute $x=n\qqq\bigpar{\tnh-\half}$ in \eqref{erika} and obtain
\begin{equation}
  \label{q2}
n\qqqw\bigpar{\xnnh-\xnh}=n\qqqw\bigpar{\xntnh-\xnh}\pto0.
\end{equation}
Furthermore, by \eqref{e1} and \eqref{e2},
\begin{equation}
  \label{q3}
\E\xnnh=\tfrac14n+O(1)
=\E\xnh+O(1).
\end{equation}
Hence, it is enough to prove \refT{T2} with $\xnnh$ replaced by
$\xnh$; we thus study
\begin{equation}\label{q4}
  \mn\=
\xxn-\xnh=
\max_{t\in \oi}\bigpar{\xnt-\xn(\thalf)}.
\end{equation}

We would like to take the supremum over all real $x$ in \eqref{erika},
but that is not allowed without further arguments since the supremum
is not a continuous functional on $D\oooo$ (the topology is too
weak). We therefore fix a large $A>0$ and study the following five
intervals separately (assuming $n>(4A)^3$): 
\begin{align*}
  I_{-2}&\=[0,\tfrac14],
\\
  I_{-1}&\=[\tfrac14,\thalf-An\qqqw],
\\
  I_{0}&\=[\thalf-An\qqqw,\thalf+An\qqqw],
\\
  I_{1}&\=[\thalf+An\qqqw,\tfrac34],
\\
  I_{2}&\=[\tfrac34,1].
\end{align*}
We denote further
\begin{equation*}
  \mnj\=\max_{t\in I_j}\bigpar{\xnt-\xn(\thalf)}_+
\end{equation*}
and have thus, since $\mn\ge0$,
\begin{equation}\label{mm}
  \mn=\max_{-2\le j\le2}\mnx{j}
\le\sum_{j=-2}^2\mnj.	
\end{equation}
On $I_0$ we use \eqref{erika}. Since the maximum is a continuous
functional on $D(I)$ for any compact interval $I$, we obtain from
\eqref{erika} on $D[-A,A]$ immediately
\begin{equation}\label{magnus}
  n\qqqw\mnx0\dto \YX_A\=\max_{|x|\le A}\bigpar{2\qqw B(x)-x^2}.
\end{equation}
Furthermore, it follows from \eqref{cecilib}, \refT{TSx} and \eqref{w3}
that
\begin{equation*}
  \E\bigpar{n\qqqw\mnx0}^2 \le C(A),
\end{equation*}
for some constant $C(A)$ depending on $A$ but not on $n$. Hence the random
variables $n\qqqw\mnx0$ are uniformly integrable, and 
\eqref{magnus} implies,
see \eg{} \cite[Theorems 5.4.2 and 5.5.9]{Gut}, 
\begin{equation}\label{k2+}
  \E\bigpar{n\qqqw\mnx0} \to \E \YX_A.
\end{equation}

On $I_{\pm2}$ we have by \eqref{cecilia}, with $\frac14\le |x|\le\half$,
\begin{equation*}
\xn\bigpar{\thalf+x}-\xn\bigpar{\thalf}
\le -n\bigpar{\tfrac14}^2+\snx1+2\snx2+6.
\end{equation*}
We use the elementary inequality, for $a>0$ and $b\in\bbR$,
\begin{equation}
  \label{jesper}
-a+b=-\frac{(a-b/2)^2}{a}+\frac{b^2}{4a}\le\frac{b^2}{4a},
\end{equation}
and obtain
\begin{equation*}
  \mnx{\pm2} \le \frac{(\snx1+2\snx2)^2}{4n/16}+6
\le 8\frac{(\snx1)^2}{n}
+32 \frac{(\snx2)^2}{n}+6
\end{equation*}
and thus by \refT{TSx}
\begin{equation}\label{k2}
\E  \mnx{\pm2}=O(1).
\end{equation}

For $I_{1}$ we define
\begin{equation}\label{k3}
  U_n(x)\=-\tfrac14\bigpar{\hsn2(\thalf+x)-\hsn2(\thalf)};
\end{equation}
this is a martingale on $[0,1/2)$.
For $0\le x\le\frac14$, we have 
\begin{equation*}
  \tfrac14\hsn2\bigpar{\thalf+x}
=
(1-2x)\qww \sn2\bigpar{\thalf+x}
= \bigpar{1+O(x)}\sn2\bigpar{\thalf+x},
\end{equation*}
and thus, using \eqref{cecilia} and \eqref{jesper}, for some constants
$C_1,C_2,\dots$, 
\begin{equation}
\label{per}
  \begin{split}
\xn\bigpar{\thalf+x}-&\xn\bigpar{\thalf}
\\&
=
-nx^2-2x\sn1\bigpar{\thalf+x}+U_n(x)
+O(x)\sn2\bigpar{\thalf+x}
+O(1)
\\
&\le
U_n(x)-\thalf nx^2+\snx1+\CC \snx2-\thalf nx^2+O(1)
\\
&\le
U_n(x)-\thalf nx^2+\CC\frac{(\snx1)^2}{n}+\CC\frac{(\snx2)^2}{n}
+O(1).
  \end{split}
\raisetag\baselineskip
\end{equation}
By \eqref{x1}, \eqref{k3} and \eqref{be22} we further have, for $0\le x\le1/4$,
\begin{equation*}
  \begin{split}
\Var\bigpar{U_n(x)}	
&=\E\qv{U_n}_x
=\tfrac1{16}\E\Bigpar{\qv{\hsn2}_{\half+x}-\qv{\hsn2}_{\half}}
\\
&=\frac{2n-2}{16}\int_{\half}^{\half+x}\frac{s}{(1-s)^3}\dd s
\le\CC nx.
  \end{split}
\end{equation*}
Hence, \refL{Lbound} yields, for $0<x_0\le\tfrac14$,
\begin{equation*}
  \E\Bigpar{\max_{x_0\le x\le 1/4}\bigpar{ U_n(x)-\thalf nx^2}_+}
\le\frac{\CC n}{n x_0}
=\frac{\CCx }{x_0}.
\end{equation*}
Taking $x_0=An\qqqw$
we thus obtain from \eqref{per},
using \refT{TSx} again, 
\begin{equation}\label{k4}
  \E\mnx1\le
\frac{\CCx }{An\qqqw}+O(1)
=
\frac{\CCx }{A}n\qqq+O(1).
\end{equation}
We obtain the same estimate for $\mnx{-1}$ by the time-reversal
$t\mapsto 1-t$ and \eqref{time}.

By \eqref{mm} and the estimates \eqref{k2} for $\mnx{\pm2}$ and
\eqref{k4} for $\mnx{\pm1}$ we find
\begin{equation*}
  \E\lrabs{\mn-\mnx0}
\le\E\mnx{-2}+\E\mnx{-1}+\E\mnx{1}+\E\mnx{2}
\le \CC+\CC n\qqq/A,
\end{equation*}
and thus
\begin{equation}\label{kk}
\limsup_\ntoo \E\bigabs{n\qqqw\mn-n\qqqw\mnx0}
\le\CCx/A.
\end{equation}

Now let $A\to\infty$; then
\begin{equation}\label{kkk}
 \YX_A\to \yoo\=\max_{x\in\bbR}\bigpar{2\qqw B(x)-x^2}.  
\end{equation}
Note that, letting $x=y/2$,
with $\YX$ as in the statement of the theorem,
\begin{equation}\label{ky}
\yoo
=\max_{y\in\bbR}\bigpar{2\qqw B(y/2)-(y/2)^2}
\eqd\max_{y\in\bbR}\bigpar{2\qw B(y)-\tfrac14 y^2}
=\thalf \YX.
\end{equation}
It follows from \eqref{kk} and \eqref{kkk} that we may let
$A\to\infty$ in \eqref{magnus} and  obtain
\begin{equation}\label{k6}
  n\qqqw\mn\dto \yoo,
\end{equation}
see \cite[Theorem 4.2]{Bill} 
(we may change the notation and denote $\mnx0$ by $\mnxx A$ for $n>(4A)^3$;
for smaller $n$ we simply let $\mnxx A=0$). 

Similarly, by \eqref{k2+} and \eqref{kk},
\begin{equation*}
\limsup_\ntoo \E\bigabs{n\qqqw\mn-\E \yoo}
\le\CCx A\qw+\bigabs{\E\yoo-\E \YX_A}.
\end{equation*}
As $A\to\infty$, $\E \YX_A\to\E\yoo$ by monotone convergence, and thus
we obtain
$\limsup_\ntoo \E\bigabs{n\qqqw\mn-\E \yoo}=0$, \ie,
\begin{equation}
  \label{k5}
n\qqqw\E\mn\to\E\yoo.
\end{equation}

The theorem follows by \eqref{k6}, \eqref{k5}, \eqref{q4}, \eqref{q2},
\eqref{q3} and \eqref{ky}.
\end{proof}

\section{Further results}\label{Sfurther}

Consider $\xnm\qi$, the number of piles with a single exam (runs of
length 1) mentioned in \refS{S:intro}.
If we for simplicity consider the cyclic case, see \refR{Rcyclic}, to
avoid edge effects (these are $O(1)$ only and do not affect the
asymptotics), we have 
\begin{equation*}
\xnm\qi=\sumkn(1-\inm(k))\inm(k+1)(1-\inm(k+2)).
\end{equation*}
After randomizing the time as in \refS{Srand}, we get
(with $\intx{k+n}=\intk$)
\begin{align}
\label{71r}
\xn\qi(t)
&=\sumkn\bigpar{1-\intk}\intx{k+1}\bigpar{1-\intx{k+2}}
\\&
=\sumkn\bigpar{1-t-\iintk}\bigpar{t+\iintx{k+1}}\bigpar{1-t-\iintx{k+2}}
\notag\\
=nt(1-t)^2&+(1-3t)(1-t)\sn1(t)-2(1-t)\sn2(t)+t\sn2'(t)+\sn3(t),
\notag
\end{align}
where we now define $\sn2$ by summing to $n$ in \eqref{sn2}, and we
introduce two new stochastic processes
\begin{align}  \label{sn2'}
  \sn2'(t)&\=\sumkn\iintk\iintx{k+2},
\\
\label{sn3}
  \sn3(t)&\=\sumkn\iintk\iintx{k+1}\iintx{k+2}.
\end{align}

The proof of \refT{TS} extends to these and yields,
in $D\oi$ and jointly with each other and \eqref{c1} and \eqref{c2}, 
\begin{align}
  \label{c2'}
n\qqw\sn2(t)&\dto Z_2'(t),
\\
  \label{c3x}
n\qqw\sn3(t)&\dto Z_3(t),
\end{align}
where $Z_2'$ and $Z_3$ are two continuous Gaussian
processes on $\oi$ with means $0$ and covariances
\begin{align}
  \E\bigpar{Z_2'(s)Z_2'(t)}&=s^2(1-t)^2,&&
 0\le s\le t\le 1,
\label{c3y}
\\
  \E\bigpar{Z_3(s)Z_3(t)}&=s^3(1-t)^3,&&
0\le s\le t\le 1.
\label{c4y}
\end{align}
Furthermore, all four processes $Z_1$, $Z_2$, $Z_2'$ and $Z_3$ are
independent. (Note that $Z_2$ and $Z_2'$ have the same distribution
but are independent.)

By the arguments in \refS{Sproofs}, which extend without any new difficulties,
this yields the following results, corresponding to our results for
$\xnm$ and $\xnt$ in Sections \refand{S:intro}{Ssn}.
We define $\xzxn:=\max_m\xnm\qi=\max_t\xn\qi(t)$, and note that
$\E\xn\qi(t)=nt(1-t)^2$ has (on $\oi$) a unique maximum at $t=1/3$.

\begin{theorem}
  \label{CXi}
As \ntoo, in $D\oi$,
\begin{multline*}
n\qqw\bigpar{\xn\qi(t)-nt(1-t)^2}
\\
\dto 
Z(t)\=(1-t)(1-3t)Z_1(t)-2(1-t)Z_2(t)+tZ_2'(t)+Z_3(t);
  \end{multline*}
$Z$ is a continuous Gaussian
process on $\oi$ with mean $\E Z(t)=0$.
\end{theorem}

\begin{theorem}
  \label{TXim}
As \ntoo, in $D\oi$,
\begin{equation*}
n\qqw\bigpar{\xnnt\qi-nt(1-t)^2}\dto Z(t)\=
-2(1-t)Z_2(t)+tZ_2'(t)+Z_3(t);
\end{equation*}
where $Z$ is a continuous Gaussian
process on $\oi$ 
with mean $\E Z(t)=0$.
\end{theorem}

We leave the explicit formulas for
(co)variances in these theorems to the reader.

\begin{theorem}
  \label{T1Xi}
As \ntoo,
\begin{equation*}
  n\qqw\bigpar{\xzxn-\tfrac4{27}n} \dto N(0,\tfrac{76}{729}),
\end{equation*}
with convergence of all moments.
In particular,
\begin{align*}
  \E\xzxn&=\frac4{27}n+o(n\qq),
\\
\Var\xzxn&=\frac{76}{729}n+o(n).
\end{align*}
\end{theorem}

\begin{theorem}
  \label{CXi2}
As \ntoo, in $D(-\infty,\infty)$,
\begin{align*}
n\qqqw\bigpar{X_n\qi(\ttt+xn\qqqw)-X_n\qi(\ttt)}\dto 
\sqrt{\frac{80}{81}}B(x)-x^2,
\end{align*}
where $B$ is a Brownian motion on
\oooo.
\end{theorem}

\begin{theorem}
  \label{T2Xi}
As \ntoo,
\begin{equation*}
  n\qqqw\bigpar{\xzxn-X_{n,\floor{n/3}}\qi}
\dto 
\gb \YX,
\end{equation*}
where the random variable $\YX$ is as in \refT{T2}
and
$\gb:
=2^{7/3}3^{-8/3}5^{2/3}=\frac{4}{27}(150)\qqq$. Furthermore,
\begin{equation*}
  \E\xzxn=
\E X_{n,\floor{n/3}}\qi+ \gb\E \YX n\qqq +o(n\qqq)
=\tfrac4{27} n+ \gb\E \YX n\qqq +o(n\qqq).
\end{equation*}
\end{theorem}

These results are easily extended to the number $\xnm\qd$ of piles
with exactly $d$ items (runs with exactly $d$ 1's) for any fixed $d$.
We may also count occurrences of any other fixed pattern, and more
generally any functional of the type
\begin{equation}
  \label{psi}
\xxx_{n,m}\=
\sumkn\psi\bigpar{\inm(k),\dots,\inm(k+\ell-1)}
\end{equation}
for some fixed $\ell\ge1$ and function $\psi:\set{0,1}^\ell\to\bbR$.
We will pursue this in some detail, leave some other details to the
reader, because the more general version illuminates the arguments
above and the structure of our method.
First, randomizing the time yields
\begin{equation}
  \label{Psi}
\xxxn(t)\=
\sumkn\Psi_k(t),
\end{equation}
where we define
$\Psi_k(t)=\psi\bigpar{\intk,\dots,\intx{k+\ell-1}}$.
We note that we will need 
more processes of the type $\sn{j}$. It turns out that it is natural to
use finite sequences of 0's and 1's to index these processes; we thus
change the notation and define a stochastic process $\snna(t)$ for
each such sequence $\ga=\ga_1\dotsm\ga_\ell$ by
\begin{equation}
  \label{snna}
\snna(t)\=\sumkn\prod_{\substack{ j\in\set{1,\dots,\ell},\\\ga_j=1}}
\iintx{k+j}.
\end{equation}
We thus now denote $\sn1(t)$, $\sn2(t)$, $\sn2'(t)$, $\sn3(t)$ by
$\snn{1}(t)$, $\snn{11}(t)$, $\snn{101}(t)$, $\snn{111}(t)$.
Initial and final 0's in $\ga$ do not affect $\snna$, so it is enough
to consider $\ga$ that begin and end with 1; let $\cA$ be the set of
all such strings $\ga$.

Let $\ell(\ga)$ denote the length of $\ga$ and $\nu(\ga)$ the
number of 1's in $\ga$, and consider only $n\ge2\ell(\ga)$. 
Then the terms in the sum in \eqref{snna} are orthogonal and we
obtain $\E\bigpar{\snna(t)}^2=n\bigpar{t(1-t)}^{\nu(\ga)}$.
Moreover, $\hsnn\ga(t):=(1-t)^{-\nu(\ga)}\snna(t)$ is a martingale on
$[0,1)$, and the proof of  \refT{TS} extends immediately and shows that,
in $D\oi$ and
jointly for all $\ga\in\cA$,
\begin{equation}\label{dsnna}
n\qqw\snna(t)\dto Z_\ga(t),	
\end{equation}
where $Z_\ga$, $\ga\in\cA$, are independent continuous Gaussian
processes with means 0 and covariances
\begin{equation}
  \E\bigpar{Z_\ga(s)Z_\ga(t)}=s^{\nu(\ga)}(1-t)^{\nu(\ga)},
\qquad 0\le s\le t\le1.
\end{equation}
Furthermore, the estimate \refT{TSx} extends to every $\snna$ (with
the implicit constant possibly depending on $\ga$).

A functional $\xxx_{n,m}$ of the type \eqref{psi} yields after
randomizing the time the functional $\xxxn(t)$ in \eqref{Psi}, which
always can be expanded as a finite sum (with orthogonal terms)
\begin{equation}\label{xxx}
  \xxxn(t)=g_0(t)n+\sum_{\ga\in\cA,\;\ell(\ga)\le\ell} g_\ga(t)\snna(t),
\end{equation}
for some polynomials $g_0(t)$ and $g_\ga(t)$, $\ga\in\cA$; this is
seen by the same argument as in \cite[Proposition 4.1]{SJ94}.
Note that, for any $n\ge2\ell$ and any $k$, 
\begin{equation}\label{g0}
 g_0(t)=\frac{\E\xxxn(t)}{n}=\E\Psi_k(t).
\end{equation}
It follows from \eqref{dsnna} and \eqref{xxx} that, 
\cf{} \refC{CX},
in $D\oi$,
\begin{equation}\label{ulla}
  n\qqw\bigpar{\xxxn(t)-ng_0(t)}
\dto Z(t)\=\sum_\ga g_\ga(t)Z_\ga(t),
\end{equation}
which is a continuous Gaussian process with mean 0 and covariance function
\begin{equation}\label{gs}
  \E\bigpar{Z(s)Z(t)} = \gs(s,t)\=
\sum_\ga g_\ga(s)g_\ga(t) (s\bmin t)^{\nu(\ga)}(1-s\bmax t)^{\nu(\ga)}.
\end{equation}
We  have moment convergence in \eqref{ulla}; moreover, the variance of
$n\qqw\xxxn(t)$ is independent of $n\ge2\ell$ and we have, for any
$k$ and $n\ge2\ell$,
\begin{equation}
  \label{gs1}
\gs(s,t)=n\qw\Cov\bigpar{\xxxn(s),\xxxn(t)}
=\sum_{j=-(\ell-1)}^{\ell-1}\Cov\bigpar{\Psi_k(s),\Psi_{k+j}(t)}.
\end{equation}

Similarly, arguing as in the proof of \refT{TXm} and
observing that the $\snn1$ terms cancel because $g_1(t)=g_0'(t)$, 
we see that, in $D\oi$,
\begin{equation}
  n\qqw\bigpar{\xxx_{n,\floor{nt}}-ng_0(t)}
\dto Z'(t)\=\sum_{\ga\neq1} g_\ga(t)Z_\ga(t),
\end{equation}
another continuous Gaussian process with mean 0.

Now suppose that $g_0(t)$ has a unique maximum on $\oi$ at an interior
point $t_0$, with $g_0''(t_0)<0$. 
Then all remaining proofs in \refS{Sproofs} extend too without
difficulties.
In particular, if we define
$\xxxxn\=\max_m\xxx_{n,m}$, we have the following generalization of \refT{T1}.
\begin{theorem}
  \label{T1x}
As \ntoo,
\begin{equation*}
  n\qqw\bigpar{\xxxxn-g(t_0)n} \dto N(0,\gss),
\end{equation*}
with convergence of all moments, with, see \eqref{gs},
\begin{equation*}
 \gss\=
\gs(t_0,t_0)=\sum_\ga g_\ga(t_0)^2\bigpar{t_0(1-t_0)}^{\nu(\ga)}.
\end{equation*}
\end{theorem}
Furthermore, \cf{} \refT{TS2},
in $D(-\infty,\infty)$ and jointly for all $\ga$,
\begin{align}
n\qqqw\bigpar{\snna(t_0+xn\qqqw)-\snna(t_0)}&\dto \gs_\ga B_\ga(x),
\end{align}
where $\gs_\ga^2\=\nu(\ga)\bigpar{t_0(1-t_0)}^{\nu(\ga)-1}$
and $B_\ga$, $\ga\in\cA$, are independent Brownian motions on \oooo.
As a consequence,
\cf{} \refC{CX2},
in $D(-\infty,\infty)$,
\begin{align}
n\qqqw\bigpar{\xxxn(t_0+xn\qqqw)-\xxxn(t_0)}\dto 
\gs\zz B(x)-\thalf|g_0''(t_0)|x^2,
\end{align}
where $B$ is a Brownian motion on \oooo{} and
\begin{equation}\label{gss}
 \gss\zz\=\sum_\ga g_\ga(t_0)^2\gss_\ga
=\sum_\ga g_\ga(t_0)^2\nu(\ga)\bigpar{t_0(1-t_0)}^{\nu(\ga)-1}. 
\end{equation}
Finally, substituting $x=(\gs\zz/|g_0''(t_0)|)\qqa y$, we see that
$\sup_{x\in\bbR} \bigpar{\gs\zz B(x)-\thalf|g_0''(t_0)|x^2}\eqd\gb\YX$,
with
\begin{equation}\label{gb}
  \gb\=(\gs\zz^4/|g_0''(t_0)|)\qqq=(\gss\zz)\qqa|g_0''(t_0)|\qqqw,
\end{equation}
and we obtain the following, where $\xxxn(t_0)$ may be replaced by
$\xxx_{n,m_0}$, where either $m_0:=\floor{t_0n}$ or $m_0$ is chosen
in \set{0,\dots,n} to maximize $\E\xxx_{n,m_0}$.
\begin{theorem}
  \label{T2x}
As \ntoo,
\begin{equation}\label{t2x}
  n\qqqw\bigpar{\xxxxn-\xxxn(t_0)}\dto 
\gb \YX,
\end{equation}
where the random variable $\YX$ is as in \refT{T2},
and $\gb$ is given by \eqref{gb} and \eqref{gss}. Furthermore,
\begin{equation*}
  \E\xxxxn=
\E\xxx(t_0)+ \gb\E \YX n\qqq +o(n\qqq)
=g_0(t_0) n+ \gb\E \YX n\qqq +o(n\qqq).
\end{equation*}
\end{theorem}

For calculation of the asymptotic variances $\gss$ and $\gss\zz$ above,
the given formulas using the coefficients $g_\ga(t_0)$ in the
decomposition \eqref{xxx} are often not very convenient.
For $\gss$, it is usually simpler to use \eqref{gs1} with $s=t=t_0$.

For $\gss\zz$ we first observe that if we take the difference of the
left derivative of $\gs(s,t)$ with respect to $s$ and the right
derivative with respect to $t$ at $(t_0,t_0)$ (thus considering $s\le
t$ only), we obtain from \eqref{gs} and \eqref{gss} easily
\begin{equation}\label{kaspar}
  \gss\zz
=
\frac{\partial}{\partial s}\gs(s,t_0)\Bigr|_{s=t_0-}
-
\frac{\partial}{\partial t}\gs(t_0,t)\Bigr|_{t=t_0+}\quad,
\end{equation}
a formula given by Daniels \cite{D89}
(in a slightly different setting).
It follows by the mean value theorem and \eqref{gs1} that, for any
fixed $n\ge2\ell$,
\begin{align}
  \gss\zz
&=
\lim_{h\downto0}\frac1h
\bigpar{\gs(t_0+h,t_0+h)-2\gs(t_0,t_0+h)+\gs(t_0,t_0)}
\\&
=
\lim_{h\downto0}\frac1{hn} \Var\bigpar{\xxxn(t_0+h)-\xxxn(t_0)}.	
\label{melchior}
  \end{align}
For fixed $n$, the probability that exactly one $\intk$ changes from 0
to 1 in the interval $[t_0,t_0+h]$ is $nh+O(h^2)$ and the probability
that more than one will change is $O(h^2)$. Hence, if $\gD_k\xxxn(t)$ is
the function of $\set{\intx{k+j}}_{1\le|j|<\ell}$ that gives the jump
in $\xxxn(t)$ (for $n\ge2\ell$) if $\intk$ is changed from 0 to 1,
keeping all other indicators fixed, then
\eqref{melchior} implies that,
for any $k$ and $n\ge2\ell$,
\begin{equation}\label{baltazar}
  \gss\zz=\Var\bigpar{\gD_k\xxxn(t_0)}.
\end{equation}

For example, for the number of runs,
\begin{equation*}
  \gD_k\xnt
=
\bigpar{1-\intx{k-1}}-\intx{k+1}
\end{equation*}
and 
$\gss\zz=\Var\bigpar{\inxx{t_0}{k-1}+\inxx{t_0}{k+1}}=2t_0(1-t_0)=\thalf$,
in accordance with \refC{CX2}.
For runs of length 1 we similarly get, from \eqref{71r},
\begin{multline*}
\gD_k\xn\qi(t)=
-\bigpar{1-\inxx{t}{k-2}}\inxx{t}{k-1}
\\
+\bigpar{1-\inxx{t}{k-1}}\bigpar{1-\inxx{t}{k}}
-\intx{k+1}\bigpar{1-\inxx{t}{k+2}}
\end{multline*}
and \eqref{baltazar} yields, in accordance with \refT{CXi2},
\begin{equation*}
\gss\zz=  \Var\bigpar{\gD_k\xn\qi(t_0)}=\frac{80}{81}.
\end{equation*}

More generally, for $\xnm\qd$ (runs of length exactly $d$), we find
from \eqref{g0},
\eqref{gs1} and \eqref{baltazar}, after straightforward calculations,
$g_0(t)=t^d(1-t)^2$, 
$t_0=d/(d+2)$,
$g_0''(t_0)=-2d^{d-1}/(d+2)^{d-1}$, and
\begin{align*}
  \gss&=\frac{4d^d}{(d+2)^{d+2}}\Bigpar{1-(d+1)\frac{4d^d}{(d+2)^{d+2}}},
\\
  \gss\zz&=8\frac{d^d}{(d+2)^{d+1}}\Bigpar{1+\frac{d^d}{(d+2)^{d+1}}},
\\
\gb
&=
\Bigpar{32\frac{d^{d+1}}{(d+2)^{d+3}}\Bigpar{1+\frac{d^d}{(d+2)^{d+1}}}^2}\qqq.
\end{align*}

\section{Priority queues and lazy hashing} \label{Spriority}

For priority queues, the $2n$ events $A_i$ and $D_i$ come in random
order, with the restriction that $A_i$ comes before $D_i$ for each
$i$. Since only the order of the events matters, we may randomize the
times as in \refS{Srand} and assume that the times $A_i$ and $D_i$,
$i=1,\dots,n$, are independent random variables 
uniformly distributed on $(0,1)$, conditioned on $A_i<D_i$ for all $i$.
For two independent random variables $T,\tT\sim U(0,1)$, the
distribution of $(T,\tT)$ conditioned on $T<\tT$ equals the
distribution of $(T\bmin\tT,T\bmax\tT)$, and this randomization of the
times in a priority queue thus gives exactly the model for lazy
hashing defined in \refS{S:intro}, as found by Kenyon and Vitter \cite{KV91}. 
In particular, $\max_t\ynt\eqd\yyn$.

In analogy with the definitions in \refS{Srand}, we now let
\begin{align*}
 \intk&\=\ett{T_k\le t},  
&
 \tintk&\=\ett{\tT_k\le t},  
\\
\iintk&\=\intk-t, 
&
\tiintk&\=\tintk-t,
\\
  \sn1(t)&\=\sumkn\iintk+\sumkn\tiintk,
\\
  \sn2(t)&\=\sumkn\iintk\tiintx{k}.
\end{align*}
We further let $N_n(t)$ be the number of events ($A_k$ or $D_k$) up to $t$.
Then, \cf{} \eqref{sofie},
\begin{align}
  \nnt
&=
\sumkn\intk+\sumkn\tintk
=\sn1(t)+2nt,
\label{ynnt}
\\
\ynt&=\sumkn\ett{A_k\le t<D_k}
=\sumkn\bigpar{\ett{T_k\le t<\tT_k}+\ett{\tT_k\le t<T_k}}
\notag
\\
&=\sumkn\Bigpar{\intk\bigpar{1-\tintk}+\tintk\bigpar{1-\intk}}
\notag
\\
&=2nt(1-t)+(1-2t)\sn1(t)-2\sn2(t).
\label{ynt}
\end{align}
We introduce martingales $\hsn1$ and $\hsn2$ as above by \eqref{ems1}
and \eqref{ems2}.

All proofs in \refS{Sproofs} now go through with no or minor changes;
the main differences are that 
\eqref{ynnt} and \eqref{ynt} contain some
factors 2 not appearing in \eqref{nnt} and \eqref{sofie}
and that
there will be a factor 2 on the \rhs{} of
\eqref{be11}; thus \refT{TS} holds with the difference that \eqref{c3}
is replaced by
\begin{equation*}
  \E\bigpar{Z_1(s)Z_1(t)}=2s(1-t),
\qquad
 0\le s\le t\le 1;
\end{equation*}
similarly, \eqref{TS2} holds with $B_1(x)$ replaced by $2\qq B_1(x)$
in \eqref{w1}. 
This yields the following results, corresponding to our results for
$\xnm$ and $\xnt$ in Sections \refand{S:intro}{Ssn}.

\begin{theorem}
  \label{CY}
As \ntoo, in $D\oi$,
\begin{equation*}
n\qqw\bigpar{\ynt-2nt(1-t)}\dto Z(t)\=(1-2t)Z_1(t)-2Z_2(t),
\end{equation*}
where $Z$ is a continuous Gaussian
process on $\oi$ with mean $\E Z(t)=0$ and covariances,
for $ 0\le s\le t\le 1$,
\begin{align*}
  \E\bigpar{Z(s)Z(t)}
&=2s(1-2s)(1-t)(1-2t)+4s^2(1-t)^2
\\
&=2s(1-t)-4s(1-s)t(1-t).
\end{align*}
\end{theorem}

\begin{theorem}
  \label{TYm}
As \ntoo, in $D\oi$,
\begin{equation*}
n\qqw\bigpar{\ynnt-2nt(1-t)}\dto Z(t)\=-2Z_2(t),
\end{equation*}
where $Z$ is a continuous Gaussian
process on $\oi$ 
with mean $\E Z(t)=0$ and covariances
\begin{align*}
  \E\bigpar{Z(s)Z(t)}&=4s^2(1-t)^2,&&
0\le s\le t\le 1.
\end{align*}
\end{theorem}

\begin{theorem}
  \label{CY2}
As \ntoo, in $D(-\infty,\infty)$,
\begin{align*}
n\qqqw\bigpar{Y_n(\thalf+xn\qqqw)-Y_n(\thalf)}\dto 2\qq B(x)-2x^2,
\end{align*}
where $B$ is a Brownian motion on
\oooo.
\end{theorem}

\begin{theorem}
  \label{T1Y}
As \ntoo,
\begin{equation*}
  n\qqw\bigpar{\yyn-n/2} \dto N(0,1/4),
\end{equation*}
with convergence of all moments.
In particular,
\begin{align*}
  \E\yyn&=n/2+o(n\qq),
\\
\Var\yyn&=n/4+o(n).
\end{align*}
\end{theorem}

\begin{theorem}
  \label{T2Y}
As \ntoo,
\begin{equation*}
  n\qqqw\bigpar{\yyn-\ynx n}\dto 
 \YX,
\end{equation*}
where the random variable $\YX$ is as in \refT{T2},
and
\begin{equation*}
  \E\yyn=
\E\ynx{n}+ \E \YX n\qqq +o(n\qqq)
=\tfrac12 n+ \E \YX n\qqq +o(n\qqq).
\end{equation*}
\end{theorem}

\refT{CY} is given by Louchard \cite{L88}, 
\refT{TYm} by Louchard \cite{L87} (with a deterministic change of
time, making the problem equivalent to a queueing problem), 
and Theorems \refand{T1Y}{T2Y} by 
Louchard, Kenyon and Schott \cite{LKS97} (with different proofs).

Note that in the proof of \refT{TYm} the terms with $\sn1$ cancel, as
discussed for 
\refT{TXm} above. In both theorems the limit is thus given by $Z_2$,
which
explains why we obtain exactly the same covariances in the two
theorems, except for a normalization factor. (Unlike \refC{CX} and
\refT{CY}, where the variances of the limits are $t(1-t)(1-3t+3t^2)$
and $2t(1-t)(1-2t+2t^2)$.)

\begin{ack}
  I thank G\"oran H\"ogn\"as for telling me about this problem and
  providing me with the reference \cite{afH}.
I also thank Anders Martin-L\"of and Guy Louchard for interesting discussions.
\end{ack}

\newcommand\AAP{\emph{Adv. Appl. Probab.} }
\newcommand\JAP{\emph{J. Appl. Probab.} }
\newcommand\JAMS{\emph{J. \AMS} }
\newcommand\MAMS{\emph{Memoirs \AMS} }
\newcommand\PAMS{\emph{Proc. \AMS} }
\newcommand\TAMS{\emph{Trans. \AMS} }
\newcommand\AnnMS{\emph{Ann. Math. Statist.} }
\newcommand\AnnPr{\emph{Ann. Probab.} }
\newcommand\CPC{\emph{Combin. Probab. Comput.} }
\newcommand\JMAA{\emph{J. Math. Anal. Appl.} }
\newcommand\RSA{\emph{Random Struct. Alg.} }
\newcommand\ZW{\emph{Z. Wahrsch. Verw. Gebiete} }
\newcommand\DMTCS{\jour{Discr. Math. Theor. Comput. Sci.} }

\newcommand\AMS{Amer. Math. Soc.}
\newcommand\Springer{Springer-Verlag}
\newcommand\Wiley{Wiley}

\newcommand\vol{\textbf}
\newcommand\jour{\emph}
\newcommand\book{\emph}
\newcommand\inbook{\emph}
\def\no#1#2,{\unskip#2, no. #1,} 
\newcommand\toappear{\unskip, to appear}

\newcommand\webcite[1]{\hfil
   \penalty0\texttt{\def~{{\tiny$\sim$}}#1}\hfill\hfill}
\newcommand\webcitesvante{\webcite{http://www.math.uu.se/~svante/papers/}}
\newcommand\arxiv[1]{\webcite{http://arxiv.org/#1}}

\def\nobibitem#1\par{}

\end{document}